\documentclass[11pt]{amsart}
\usepackage[a4paper,margin=2.5cm]{geometry}
\usepackage[T1]{fontenc}
\usepackage{lmodern}
\usepackage{amssymb,mathtools}
\usepackage{graphicx,booktabs,multirow}
\usepackage{enumitem,placeins}
\usepackage{natbib}
\usepackage[hidelinks]{hyperref}

\numberwithin{equation}{section}

\theoremstyle{plain}
\newtheorem{theorem}{Theorem}[section]
\newtheorem{proposition}[theorem]{Proposition}
\newtheorem{corollary}[theorem]{Corollary}
\newtheorem{lemma}[theorem]{Lemma}
\theoremstyle{definition}
\newtheorem{assumption}[theorem]{Assumption}
\theoremstyle{remark}

\newcommand{\R}{\mathbb R}
\newcommand{\E}{\mathbb E}
\newcommand{\PP}{\mathbb P}

\newcommand{\cD}{\mathcal D}
\newcommand{\cF}{\mathcal F}
\newcommand{\cG}{\mathcal G}
\newcommand{\cP}{\mathcal P}
\newcommand{\argmin}{\mathop{\rm arg\,min}}

\begin{document}

\title{Parareal Approximation of Stochastic $L^2$-Boosting}
\author{Duc-Quang Bui}
\address{Laboratoire Manceau de Math\'ematiques, Le Mans Universit\'e, France}
\email{duc\_quang.bui@univ-lemans.fr}

\begin{abstract}
We study Parareal approximations of stochastic \(L^2\)-boosting, reusing the learner draws of the sequential algorithm in the fine propagators. We establish mean-square error bounds whose order improves with each correction, characterize the limiting distributions of normalized errors, including non-Gaussian components and the effects of coarse discretization, and derive sufficient conditions for preserving the statistical limit laws and validation-based stopping behavior of the sequential algorithm. Numerical experiments with sample-dependent coordinate and nonlinear residual-fitting learners examine the predicted convergence rates and limiting distributions, while computational benchmarks assess the efficiency of the method relative to sequential implementations.
\end{abstract}


\maketitle

\section{Introduction}

Stochastic gradient boosting randomizes the learner construction, for instance through subsampling, to reduce computational cost and improve prediction \citep{Friedman2002}. Its vanishing-learning-rate behavior has been studied through infinitesimal boosting flows for linear learners and randomized trees \citep{DombryEsstafa2024,DombryDuchamps2024}, together with their large-sample behavior and fluctuations \citep{DombryDuchampsLarge2024,DombryDuchampsFCLT2026}. Related diffusion and gradient-flow analyses provide limit theorems and inference for stochastic algorithms \citep{Latz2021,FlamandBrunel2026,LiGiessing2026}. For stochastic $L^2$-boosting, \citet{BuiEsstafaJoint2026} establish joint functional limits for sampling and learner fluctuations and identify the scales governing validation-based early stopping.

The Parareal method of \citet{LionsMadayTurinici2001} combines sequential coarse propagation with fine computations performed concurrently over time subintervals; its convergence was analyzed by \citet{GanderVandewalle2007}. It can also be coupled with other solvers, including optimized Schwarz waveform relaxation \citep{BuiJaphetMadayOmnes2022}. In mathematical finance, \citet{BalMaday2002} apply Parareal to the PDE formulation of American-option pricing, while \citet{PagesPironneauSall2018} combine it with least-squares Monte Carlo. Stochastic and probabilistic developments include randomized Parareal corrections (Pentland, Tamborrino and Sullivan \citeyear{PentlandTamborrinoSullivan2023}), Gaussian-process emulation of the fine--coarse discrepancy \citep{PentlandEtAl2023GParareal,GattiglioEtAl2026}, and parallel-in-time probabilistic ODE solvers based on filtering and smoothing \citep{BoschEtAl2024}.

Related approaches connect time integration with optimization. ParaFlow uses Parareal corrections to accelerate gradient-flow minimization \citep{BellezzaEtAl2026}, while \citet{SunXu2026} consider orthogonality-constrained gradient flows. Other methods parallelize stochastic optimization by solving triangular nonlinear systems \citep{LuEtAl2026PASO}, or exploit associative scans within Newton solvers for ODEs \citep{ParallelNewtonScan2026}.

Our fine propagators reuse the learner draws from the sequential algorithm on every Parareal correction. This coupling permits a conditional mean-square error estimate whose constant is independent of the number of subintervals. For smooth randomized learners, we establish limit distributions of the normalized errors at fixed correction levels and distinguish the effects of learner randomness and coarse Euler bias. We then derive sufficient conditions for the numerical error to be negligible at the statistical path and validation-selection scales. The numerical study examines sample-dependent coordinate and nonlinear residual-fitting learners, including their second-order error distributions; measured timings are compared with sequential algorithms that include validation.

The paper is organized as follows. Section~\ref{sec:setting} introduces the coupled algorithm. Sections~\ref{sec:bounds} and~\ref{sec:second} establish error bounds and limiting error laws. Section~\ref{sec:statistics} studies statistical accuracy and computational cost, and Section~\ref{sec:numerical} presents numerical experiments.

\section{Randomized \texorpdfstring{$L^2$}{L2}-boosting and Parareal}
\label{sec:setting}

Fix a time horizon $T>0$ and an integer $M\geq2$. Partition $[0,T]$ into $M$ equal subintervals $[T_j,T_{j+1}]$, where
\[
\Delta=\frac{T}{M}, \qquad T_j=j\Delta,\quad 0\leq j\leq M.
\]
Each subinterval contains $L\geq1$ fine time steps, with boosting learning rate
\[
\lambda=\frac{\Delta}{L}.
\]
The total number of fine steps is $N=ML$. The coarse propagator uses $J\geq1$ explicit Euler steps per subinterval, with step size $\tau=\Delta/J$.

Let $n$ denote the training-sample size and $p\geq1$ the dimension of the boosting state. We suppress the dependence on $n$ in the algorithmic notation; the relationships between the sample size and the discretization parameters will be specified in the asymptotic analysis.

Let $\cD_n$ be the sigma-field generated by the training observations, which may be dependent. At each boosting step, we average $B\geq1$ learner draws. Let $\{\xi_{\ell,r}:1\leq\ell\leq N,\ 1\leq r\leq B\}$ be independent and identically distributed conditionally on $\cD_n$. We write $\E_\xi$ for expectation over these draws, conditionally on $\cD_n$.

For a state $x\in\R^p$ and a generic learner draw $\xi$, let $\widehat F(x,\xi)\in\R^p$ be a measurable, conditionally integrable update direction. The update direction may depend on the training data. Its conditional mean is
\[
f(x)=\E_\xi\widehat F(x,\xi).
\]
For a deterministic initial state $\beta_0\in\R^p$, the sequential boosting iterates are defined by
\begin{equation}\label{eq:general-boost}
\begin{aligned}
\Theta_0&=\beta_0,\\
\Theta_{\ell+1} &=\Theta_\ell +\frac{\lambda}{B}\sum_{r=1}^{B} \widehat F(\Theta_\ell,\xi_{\ell+1,r}), \qquad 0\leq\ell<N.
\end{aligned}
\end{equation}
Let $\Theta$ denote their piecewise-linear interpolation, satisfying
\[
\Theta(\ell\lambda)=\Theta_\ell, \qquad 0\leq\ell\leq N.
\]

For each $j=0,\ldots,M-1$ and $x\in\R^p$, define the fine iterates
\[
z_0=x,\qquad z_{\ell+1} =z_\ell+\frac{\lambda}{B} \sum_{r=1}^{B} \widehat F(z_\ell,\xi_{jL+\ell+1,r}), \quad 0\leq\ell<L.
\]
The fine propagator over $[T_j,T_{j+1}]$ is
\[
\cF_j(x)=z_L.
\]
Let $\cP_j(t;x)$, $t\in[T_j,T_{j+1}]$, be the piecewise-linear interpolation of these iterates, satisfying
\[
\cP_j(T_j+\ell\lambda;x)=z_\ell, \qquad 0\leq\ell\leq L.
\]
For $0\leq u\leq\Delta$, define the partial fine propagator by
\[
\cF_{j,u}(x)=\cP_j(T_j+u;x).
\]
In particular,
\[
\cF_{j,0}(x)=x,\qquad \cF_{j,\Delta}(x)=\cF_j(x).
\]
The sequential stochastic path therefore satisfies
\[
\Theta(t) =\cP_j\bigl(t;\Theta(T_j)\bigr), \qquad t\in[T_j,T_{j+1}].
\]
The learner draws are generated once and reused in every evaluation of the fine propagators. Consequently, these propagators are fixed conditionally on the training data and learner draws.

The coarse propagator approximates the flow of the conditional mean field $f$ by the explicit Euler method. For an initial state $x$, define
\[
y_0=x,\qquad y_{k+1}=y_k+\tau f(y_k), \quad 0\leq k<J,
\]
and set
\[
\cG(x)=y_J.
\]
Conditionally on $\cD_n$, the coarse propagator is deterministic. The fine propagator reproduces the prescribed stochastic boosting updates, whereas the coarse propagator approximates the corresponding mean dynamics.

For $K\geq0$, let $U_j^{(K)}$ denote the Parareal approximation at $T_j$. The initial coarse prediction is
\[
U_0^{(0)}=\beta_0,\qquad U_{j+1}^{(0)}=\cG(U_j^{(0)}), \quad 0\leq j<M.
\]
Successive corrections are defined by
\begin{equation}\label{eq:parareal}
\begin{aligned}
U_0^{(K+1)}&=\beta_0,\\
U_{j+1}^{(K+1)} &=\cG(U_j^{(K+1)}) +\cF_j(U_j^{(K)}) -\cG(U_j^{(K)}), \quad 0\leq j<M.
\end{aligned}
\end{equation}
At each correction, the fine propagations $\cF_j(U_j^{(K)})$ can be computed concurrently over the time subintervals, followed by a sequential coarse sweep.

A continuous Parareal approximation is obtained by correcting the fine interpolation on each subinterval:
\begin{equation}\label{eq:bridge}
\Theta^{(K,M)}(T_j+u) =\cF_{j,u}(U_j^{(K)})+\frac{u}{\Delta} \left\{ U_{j+1}^{(K)} -\cF_j(U_j^{(K)}) \right\}, \qquad 0\leq u\leq\Delta.
\end{equation}
The correction matches the Parareal endpoints, so that
\[
\Theta^{(K,M)}(T_j)=U_j^{(K)}, \qquad 0\leq j\leq M,
\]
and the reconstructed path is continuous.

\noindent \textbf{Example. Randomized coordinate boosting.} For training observations $(X_i,Y_i)\in\R^p\times\R$, $i=1,\ldots,n$, define
\[
H_n=\frac1n\sum_{i=1}^n X_iX_i^\top, \qquad b_n=\frac1n\sum_{i=1}^n X_iY_i.
\]
Choose probabilities $\pi_1,\ldots,\pi_p>0$ summing to one. Let $\iota$ be independent of the training sample, with $\PP(\iota=i)=\pi_i$, and let $e_i$ be the $i$th coordinate vector. Define
\[
D(\iota)=\pi_\iota^{-1}e_\iota e_\iota^\top, \qquad \widehat F(x,\iota)=D(\iota)(b_n-H_nx).
\]
Since $\E_\xi D(\iota)=I_p$, where $I_p$ denotes the identity matrix, the conditional mean direction is
\[
f(x)=b_n-H_nx.
\]

Let $\iota_{\ell,r}$ be independent copies of $\iota$ and set
\[
\widehat D_\ell =\frac1B\sum_{r=1}^B D(\iota_{\ell,r}), \qquad 1\leq\ell\leq N.
\]
Then \eqref{eq:general-boost} becomes
\begin{equation}\label{eq:boost}
\Theta_{\ell+1} =\Theta_\ell+\lambda\widehat D_{\ell+1} (b_n-H_n\Theta_\ell), \qquad 0\leq\ell<N.
\end{equation}
The corresponding coarse Euler updates are
\[
y_{k+1}=y_k+\tau(b_n-H_ny_k), \qquad 0\leq k<J.
\]

\section{Error bounds}\label{sec:bounds}
We first record the finite-step exactness of Parareal. Each correction makes the approximation exact on one additional subinterval. Thus, after $M$ corrections, the algorithm reproduces the entire sequential trajectory. The error bounds then quantify the accuracy obtained with fewer corrections.

\begin{proposition}\label{prop:finite-exactness}
For every integer $K\ge0$, the Parareal endpoints satisfy
\[
U_{j}^{(K)}=\Theta(T_j), \qquad 0\le j\le\min\{K,M\}.
\]
Consequently, for every $K\ge M$,
\[
\Theta^{(K,M)}(t)=\Theta(t), \qquad 0\le t\le T.
\]
These identities hold pathwise.
\end{proposition}

\begin{proof}
All correction levels start from $\beta_0$. Suppose level $K$ is exact through node $\min\{K,M\}$. Starting from node zero, an induction over the subintervals gives
\[
\begin{aligned}
U_{j+1}^{(K+1)} &= \mathcal G(\Theta(T_j)) +\mathcal F_{j}(\Theta(T_j)) -\mathcal G(\Theta(T_j))\\
&=\Theta(T_{j+1}), \qquad 0\le j\le\min\{K,M-1\}.
\end{aligned}
\]
This proves the endpoint identities by induction on $K$. The fine maps use the same learner draws as the sequential trajectory, so their evaluations agree pathwise.

For $K\ge M$, every endpoint is exact. The endpoint correction in \eqref{eq:bridge} therefore vanishes, and the partial fine maps reproduce the sequential path on every subinterval.
\end{proof}

\subsection{General error bounds}
\label{subsec:general-bounds}

We study the difference between the Parareal reconstruction and the sequential stochastic boosting path. All expectations denoted by $\E_\xi$ are conditional on $\cD_n$.

We use Euclidean vector norms and induced operator norms for matrices. The Frobenius norm is denoted by $\|\cdot\|_F$. For a differentiable map $g$, $Dg$ denotes its derivative, and $I_p$ denotes the $p\times p$ identity matrix. For a random vector or scalar $Y$ and an exponent $q\ge1$, write
\[
\|Y\|_{L^q_\xi} = \bigl(\E_\xi\|Y\|^q\bigr)^{1/q}.
\]
For a continuous path $v:[0,T]\to\R^p$, set $\|v\|_\infty=\sup_{0\le t\le T}\|v(t)\|$.

For $j=0,\ldots,M$, define
\[
V_{j}=\Theta(T_j), \qquad \mathcal H_{j} = \sigma\!\left( \cD_n,\, \xi_{\ell,r}: 1\le\ell\le jL,\ 1\le r\le B \right).
\]
Thus $\mathcal H_{j}$ contains the training data and the learner draws in the subintervals preceding $T_j$. The draws in subinterval $j$ are independent of $\mathcal H_{j}$ conditionally on $\cD_n$. Define the reference learner-noise scale by
\[
s_n=\sqrt{\frac{\lambda}{B}}.
\]

\begin{assumption}\label{ass:nonlin-propagator}
For every deterministic $x\in\R^p$ and every $j=0,\ldots,M-1$,
\[
\E_\xi\|\cF_{j}(x)\|^2<\infty.
\]
Define
\[
\mathcal Z_{j}(x) = \cF_{j}(x)-\E_\xi\cF_{j}(x), \qquad \mathcal B_{j}(x) = \E_\xi\cF_{j}(x)-\cG(x).
\]
There exist a nonnegative deterministic sequence $d_n$ and a constant $C\ge1$, independent of $n$, $M$, $L$, $J$ and $B$, such that the following bounds hold conditionally on $\cD_n$, simultaneously for all starting states.
\begin{enumerate}[label=(\roman*)]
\item The coarse map is continuously differentiable and
\[
\sup_{x\in\R^p} \|D\cG(x)-I_p\| \le C\Delta.
\]

\item For $j=0,\ldots,M-1$ and $x,y\in\R^p$,
\begin{align}
\E_\xi \|\mathcal Z_{j}(x)-\mathcal Z_{j}(y)\|^2 &\le C\Delta s_n^2\|x-y\|^2, \tag{N1}\\
\E_\xi\|\mathcal Z_{j}(x)\|^2 &\le C\Delta s_n^2(1+\|x\|^2), \tag{N2}\\
\|\mathcal B_{j}(x)-\mathcal B_{j}(y)\| &\le C\Delta d_n\|x-y\|, \tag{N3}\\
\|\mathcal B_{j}(x)\| &\le C\Delta d_n(1+\|x\|). \tag{N4}
\end{align}

\item Almost surely, for every $j=0,\ldots,M-1$, $u\in[0,\Delta]$ and $x,y\in\R^p$,
\[
\|\cF_{j,u}(x)-\cF_{j,u}(y)\| \le C\|x-y\|.
\]

\item The sequential endpoints satisfy
\[
\E_\xi \max_{0\le j\le M}(1+\|V_{j}\|)^2 \le C.
\]
\end{enumerate}
\end{assumption}

The fine--coarse defect decomposes as
\[
\cF_{j}-\cG = \mathcal Z_{j}+\mathcal B_{j}.
\]
The centered term describes learner fluctuations, while the second term measures the deterministic subinterval mismatch conditional on the data. The coupled increment bound compares two starting states using the same learner draws. It provides the additional noise factor at each correction. The bound on the partial fine maps transfers endpoint estimates to the reconstructed continuous path.

We first control the propagation of a centered perturbation together with a bounded drift.

\begin{lemma}\label{lem:nonlin-predictable}
Assume part (i) of Assumption~\ref{ass:nonlin-propagator}. Fix $n$. Let $u_j$, $j=0,\ldots,M$, and $v_j$, $j=0,\ldots,M-1$, be $\mathcal H_{j}$-measurable vectors, with $u_0=0$, satisfying
\[
u_{j+1} = \cG(v_j+u_j)-\cG(v_j)+z_j+b_j.
\]
For $j=0,\ldots,M-1$, let $a_j\ge0$ be $\mathcal H_{j}$-measurable, with $\|\max_j a_j\|_{L^2_\xi}<\infty$. Assume that $z_j$ and $b_j$ are $\mathcal H_{j+1}$-measurable and
\[
\E(z_j\mid\mathcal H_{j})=0, \qquad \E(\|z_j\|^2\mid\mathcal H_{j}) \le C\Delta s_n^2a_j^2, \qquad \|b_j\|\le C\Delta d_na_j.
\]
Then
\[
\left\|\max_{0\le j\le M}\|u_j\|\right\|_{L^2_\xi} \le c_T(s_n+d_n) \left\|\max_{0\le j<M}a_j\right\|_{L^2_\xi},
\]
where $c_T$ depends only on $T$ and $C$.
\end{lemma}

\begin{proof}
The mean-value theorem gives
\[
\cG(v_j+u_j)-\cG(v_j)=W_ju_j, \qquad W_j=\int_0^1D\cG(v_j+au_j)\,da,
\]
with $\|W_j-I_p\|\le C\Delta$. Summing the recursion yields
\[
u_i = \sum_{j<i}(W_j-I_p)u_j +\sum_{j<i}b_j +\sum_{j<i}z_j.
\]
Set
\[
S=\max_{0\le i\le M} \left\|\sum_{j<i}z_j\right\|, \qquad \bar a=\max_{0\le j<M}a_j.
\]
Since $M\Delta=T$,
\[
\|u_i\| \le C\Delta\sum_{j<i}\|u_j\|+CTd_n\bar a+S.
\]
Discrete Gronwall therefore gives
\[
\max_{0\le i\le M}\|u_i\| \le e^{CT}(CTd_n\bar a+S).
\]
The partial sums of $z_j$ form an $(\mathcal H_{i})$-martingale. Conditional Doob's inequality and orthogonality imply
\[
\E_\xi S^2 \le 4\sum_{j<M}\E_\xi\|z_j\|^2 \le 4CTs_n^2\E_\xi\bar a^2.
\]
Taking conditional $L^2$ norms proves the assertion, for example with
\[
c_T=e^{CT}\bigl(2\sqrt{CT}+CT\bigr).\qedhere
\]
\end{proof}

Applying this estimate to the Parareal errors gives the following bound.

\begin{theorem}\label{thm:nonlinear-uniform}
Under Assumption~\ref{ass:nonlin-propagator}, there exists $C_T\ge1$ such that, for every $K\ge0$,
\begin{equation}\label{eq:nonlin-bound}
\left( \E_\xi \max_{0\le j\le M} \|U_{j}^{(K)}-V_{j}\|^2 \right)^{1/2} \le C_T^{K+1}(s_n+d_n)^{K+1}.
\end{equation}
Moreover,
\begin{equation}\label{eq:nonlin-full}
\left( \E_\xi \|\Theta^{(K,M)}-\Theta\|_\infty^2 \right)^{1/2} \le (1+2C)C_T^{K+1}(s_n+d_n)^{K+1}.
\end{equation}
The constant $C_T$ depends only on $T$ and the constants of Assumption~\ref{ass:nonlin-propagator}. It is independent of $n$, $M$, $L$, $J$, $B$ and $K$. For $K\ge M$, the reconstructed path equals $\Theta$.
\end{theorem}

\begin{proof}
Define
\[
e_{j}^{(K)}=U_{j}^{(K)}-V_{j}, \qquad A_K= \left( \E_\xi\max_j\|e_{j}^{(K)}\|^2 \right)^{1/2}.
\]
Subtracting the sequential fine recursion from \eqref{eq:parareal} gives
\begin{equation}\label{eq:nonlin-exact}
\begin{aligned}
e_{j+1}^{(K+1)} ={}& \cG(V_{j}+e_{j}^{(K+1)}) -\cG(V_{j})\\
&+ \mathcal Z_{j}(V_{j}+e_{j}^{(K)}) -\mathcal Z_{j}(V_{j})\\
&+ \mathcal B_{j}(V_{j}+e_{j}^{(K)}) -\mathcal B_{j}(V_{j}).
\end{aligned}
\end{equation}

Induction in the Parareal recursion shows that $U_{j}^{(K)}$, $V_{j}$ and $e_{j}^{(K)}$ are $\mathcal H_{j}$-measurable. The current subinterval draws are independent of this sigma-field conditionally on the data. Hence the centered difference in \eqref{eq:nonlin-exact} has conditional mean zero, and its conditional second moment is bounded by $C\Delta s_n^2\|e_{j}^{(K)}\|^2$. The bias difference is bounded by $C\Delta d_n\|e_{j}^{(K)}\|$. Lemma~\ref{lem:nonlin-predictable} yields
\[
A_{K+1}\le c_T(s_n+d_n)A_K.
\]
The reuse of learner draws ensures that the coupled increment bound applies to two evaluations of the same random propagator.

For the coarse predictor,
\[
e_{j+1}^{(0)} = \cG(V_{j}+e_{j}^{(0)}) -\cG(V_{j}) -\mathcal Z_{j}(V_{j}) -\mathcal B_{j}(V_{j}).
\]
The same lemma, with $a_j=1+\|V_{j}\|$, gives
\[
A_0\le c_T\sqrt C\,(s_n+d_n).
\]
It follows that
\[
A_K\le \sqrt C\,c_T^{K+1}(s_n+d_n)^{K+1}.
\]
Choosing $C_T\ge\max\{1,c_T,c_T\sqrt C\}$ proves \eqref{eq:nonlin-bound}.

On subinterval $j$, subtract the sequential fine path from \eqref{eq:bridge}. The Lipschitz bounds on the partial fine maps give
\[
\|\Theta^{(K,M)}(T_j+u)-\Theta(T_j+u)\| \le 2C\|e_{j}^{(K)}\|+\|e_{j+1}^{(K)}\|.
\]
Taking the supremum over $t\in[0,T]$ proves \eqref{eq:nonlin-full}.

Finally, induction gives $e_{j}^{(K)}=0$ for $j\le K$. All endpoints agree for $K\ge M$, and the reconstruction then reproduces the sequential fine path.
\end{proof}

The same argument gives a weaker improvement when the coupled defect is only H\"older continuous.

\begin{theorem}\label{thm:holder}
Let $\alpha\in(0,1]$. Retain the conditional square-integrability requirement, parts (i) and (iv), and bounds (N2) and (N4) of Assumption~\ref{ass:nonlin-propagator}. Replace (N1) and (N3) by
\begin{equation}\label{eq:holder-cond}
\begin{aligned}
\E_\xi \|\mathcal Z_{j}(x)-\mathcal Z_{j}(y)\|^2 &\le C\Delta s_n^2\|x-y\|^{2\alpha},\\
\|\mathcal B_{j}(x)-\mathcal B_{j}(y)\| &\le C\Delta d_n\|x-y\|^\alpha.
\end{aligned}
\end{equation}
For each fixed $K\ge0$,
\begin{equation}\label{eq:holder-rate}
\left( \E_\xi\max_{0\le j\le M} \|U_{j}^{(K)}-V_{j}\|^2 \right)^{1/2} \le C_{T,K,\alpha}(s_n+d_n)^{\nu_K}, \qquad \nu_K=\sum_{r=0}^{K}\alpha^r.
\end{equation}
The constant is independent of $n$, $M$, $L$, $J$ and $B$.
\end{theorem}

\begin{proof}
Fix $n$ and define
\[
A_K= \left( \E_\xi\max_j\|e_{j}^{(K)}\|^2 \right)^{1/2}.
\]
The predictor estimate remains $A_0\le C_T(s_n+d_n)$. Apply Lemma~\ref{lem:nonlin-predictable} to \eqref{eq:nonlin-exact} with $a_j=\|e_{j}^{(K)}\|^\alpha$. Jensen's inequality gives
\[
\begin{aligned}
A_{K+1} &\le c_T(s_n+d_n) \left( \E_\xi\max_j\|e_{j}^{(K)}\|^{2\alpha} \right)^{1/2}\\
&\le c_T(s_n+d_n)A_K^\alpha.
\end{aligned}
\]
Iteration proves the result, with $\nu_0=1$ and $\nu_{K+1}=1+\alpha\nu_K$.
\end{proof}

\subsection{Smooth learners}
\label{subsec:smooth-learners}

We next verify the propagator assumptions from regularity properties of the learner updates. For an integer $k\ge0$ and a field $g$ with derivatives up to order $k$, write
\[
\|g\|_{C^k} = \sum_{a=0}^{k}\sup_{x\in\R^p}\|D^ag(x)\|, \qquad D^0g=g.
\]
For a matrix-valued field $\sigma$, $D\sigma(x)[v]$ denotes its directional derivative along $v\in\R^p$.

\begin{assumption}\label{ass:primitive}
There is a fixed integer $m\ge1$ such that the learner draws take values in $\R^m$ and
\[
\widehat F(x,\xi)=f(x)+\sigma(x)\xi.
\]
The mean field $f:\R^p\to\R^p$ and the $\cD_n$-measurable field $\sigma:\R^p\to\R^{p\times m}$ are $C^2$. There exist deterministic constants $C_f,L_\xi<\infty$, independent of $n$, such that
\[
\|f\|_{C^2}+\|\sigma\|_{C^2}\le C_f
\]
and, conditionally on $\cD_n$,
\[
\E_\xi\xi_{\ell,r}=0, \qquad \|\xi_{\ell,r}\|\le L_\xi \quad\text{almost surely}.
\]
\end{assumption}

For $\ell=1,\ldots,N$, define the batch average
\[
\bar\xi_{\ell} = \frac1{B}\sum_{r=1}^{B}\xi_{\ell,r}.
\]
The fine recursion becomes
\begin{equation}\label{eq:nonlin-primitive}
\Theta_{\ell+1} = \Theta_{\ell} +\lambda \left\{ f(\Theta_{\ell}) +\sigma(\Theta_{\ell})\bar\xi_{\ell+1} \right\}.
\end{equation}
Define the combined scale
\[
\eta_n=s_n+\lambda+\tau.
\]

\begin{proposition}\label{prop:primitive}
Under Assumption~\ref{ass:primitive}, the fine maps and the multistep Euler coarse map of Section~\ref{sec:setting} satisfy Assumption~\ref{ass:nonlin-propagator} with
\[
d_n=\lambda+\tau+s_n\sqrt{\Delta}.
\]
There exists $C_T\ge1$ such that, for every $K\ge0$,
\[
\left( \E_\xi \|\Theta^{(K,M)}-\Theta\|_\infty^2 \right)^{1/2} \le C_T^{K+2}\eta_n^{K+1}.
\]
The constant depends on $T$, $p$, $m$, $C_f$, $L_\xi$ and $\|\beta_0\|$, but not on $n$, $M$, $L$, $J$, $B$ or $K$.
\end{proposition}

\begin{proof}
Fix $n$, a subinterval $j$ and a deterministic starting state $x$. Within this proof, put $\zeta_\ell=\bar\xi_{jL+\ell}$ for $\ell=1,\ldots,L$. Conditional independence and centering give
\[
\E_\xi\zeta_\ell=0, \qquad \E_\xi\|\zeta_\ell\|^2\le \frac{L_\xi^2}{B}, \qquad \|\zeta_\ell\|\le L_\xi.
\]
Let $X_\ell(x)$ and $\bar X_\ell(x)$ denote the local fine trajectory and the noise-free Euler trajectory. They satisfy
\[
\begin{aligned}
X_0(x)&=\bar X_0(x)=x,\\
X_{\ell+1} &= X_\ell+\lambda \{f(X_\ell)+\sigma(X_\ell)\zeta_{\ell+1}\},\\
\bar X_{\ell+1} &= \bar X_\ell+\lambda f(\bar X_\ell), \qquad 0\le\ell<L.
\end{aligned}
\]
Write
\[
P_\ell=D_xX_\ell(x), \qquad \bar P_\ell=D_x\bar X_\ell(x).
\]
Both Jacobians start from $I_p$. Constants denoted by $C$ in this proof may depend on $T$, $p$, $m$, $C_f$ and $L_\xi$, but not on the discretization levels or $B$.

For $v\in\R^p$, the fine Jacobian recursion is
\[
\begin{aligned}
P_{\ell+1}v ={}&P_\ell v +\lambda Df(X_\ell)P_\ell v\\
&+\lambda D\sigma(X_\ell)[P_\ell v]\zeta_{\ell+1}.
\end{aligned}
\]
Each fine update therefore has derivative norm bounded by $1+C\lambda$. It follows that
\[
\|P_\ell\|+\|\bar P_\ell\|\le C, \qquad 0\le\ell\le L.
\]
The same Lipschitz bound holds for the linearly interpolated partial fine maps.

The fine increments have norm at most $C\lambda$. Hence $\|\cF_{j}(x)\|\le\|x\|+C\Delta$, which proves conditional square integrability at every deterministic starting state. The complete sequential path is bounded by $\|\beta_0\|+CT$. This verifies parts (iii) and (iv) of the propagator assumption.

For the coarse trajectory $y_k(x)$, define $R_k=D_xy_k(x)$. Then
\[
R_0=I_p, \qquad R_{k+1}=R_k+\tau Df(y_k)R_k.
\]
Thus $\|R_k\|\le e^{C_f\Delta}$ and
\[
D\cG(x)-I_p = \tau\sum_{k=0}^{J-1}Df(y_k)R_k.
\]
Since $J\tau=\Delta$, this proves part (i).

To estimate the centered fluctuations, define
\[
E_\ell=X_\ell-\bar X_\ell, \qquad \Gamma_\ell=P_\ell-\bar P_\ell.
\]
Both differences start from zero. The state difference satisfies
\[
E_{\ell+1} = E_\ell+\lambda \{f(X_\ell)-f(\bar X_\ell)\} +\lambda\sigma(X_\ell)\zeta_{\ell+1}.
\]
The last term is centered conditionally on the preceding local draws. Consequently,
\[
\E_\xi\|E_{\ell+1}\|^2 \le (1+C\lambda)\E_\xi\|E_\ell\|^2 +C\frac{\lambda^2}{B}.
\]
Discrete Gronwall gives
\begin{equation}\label{eq:prim-E}
\E_\xi\|E_\ell\|^2 \le C\frac{\ell\lambda^2}{B} = C(\ell\lambda)s_n^2.
\end{equation}

The drift in the Jacobian-difference recursion is
\[
Df(X_\ell)P_\ell -Df(\bar X_\ell)\bar P_\ell.
\]
Its Frobenius norm is bounded by $C(\|\Gamma_\ell\|_F+\|E_\ell\|)$. The Jacobian-noise increment is conditionally centered and has conditional second moment bounded by $C\lambda^2/B$. Therefore,
\[
\E_\xi\|\Gamma_{\ell+1}\|_F^2 \le (1+C\lambda)\E_\xi\|\Gamma_\ell\|_F^2 +C\lambda\E_\xi\|E_\ell\|^2 +C\frac{\lambda^2}{B}.
\]
Using \eqref{eq:prim-E} and $\ell\lambda\le T$, another application of discrete Gronwall yields
\begin{equation}\label{eq:prim-Gamma}
\E_\xi\|\Gamma_\ell\|_F^2 \le C(\ell\lambda)s_n^2.
\end{equation}

Since $\bar X_{L}(x)$ is deterministic conditionally on the data,
\[
\E_\xi\|\mathcal Z_{j}(x)\|^2 \le \E_\xi\|X_{L}(x)-\bar X_{L}(x)\|^2 \le C\Delta s_n^2.
\]
This proves (N2). The pathwise Jacobian bound permits differentiation under conditional expectation. For $x,y\in\R^p$, set $x_a=y+a(x-y)$ for $a\in[0,1]$. Then
\[
\mathcal Z_{j}(x)-\mathcal Z_{j}(y) = \int_0^1 \{P_{L}(x_a)-\E_\xi P_{L}(x_a)\}(x-y)\,da.
\]
The variance-minimizing property of conditional expectation and \eqref{eq:prim-Gamma} give
\[
\E_\xi \|P_{L}(x_a)-\E_\xi P_{L}(x_a)\|_F^2 \le C\Delta s_n^2.
\]
The Cauchy--Schwarz inequality now proves (N1).

It remains to estimate the deterministic mismatch. For $h>0$ with $\Delta/h\in\mathbb N$, let $\mathcal E_{h}$ denote deterministic Euler propagation over time $\Delta$ for the field $f$, using step $h$. In particular,
\[
\mathcal E_{\lambda}(x)=\bar X_{L}(x), \qquad \mathcal E_{\tau}=\cG.
\]
Taking expectations in the recursions for $E_\ell$ and $\Gamma_\ell$ removes the centered noise terms. The drift bounds and \eqref{eq:prim-E}--\eqref{eq:prim-Gamma} imply
\[
\begin{aligned}
\|\E_\xi E_{L}(x)\| +\|\E_\xi\Gamma_{L}(x)\|_F & \le C\lambda\sum_{\ell=0}^{L-1} \left\{ (\E_\xi\|E_\ell(x)\|^2)^{1/2} + (\E_\xi\|\Gamma_\ell(x)\|_F^2)^{1/2} \right\}\\
&\le Cs_n\lambda \sum_{\ell=0}^{L-1}\sqrt{\ell\lambda} \le C\Delta^{3/2}s_n.
\end{aligned}
\]
The estimates are uniform in $x$, so
\[
\|\E_\xi\cF_{j} -\mathcal E_{\lambda}\|_{C^1} \le C\Delta^{3/2}s_n.
\]

Let $\phi_{t}(x)$ denote the exact solution at time $t$ of $\dot x=f(x)$, started from $x$. The global $C^2$ bounds ensure existence of this flow and bounded state derivatives on $[0,T]$. The local Euler errors for the state and its Jacobian are bounded by $Ch^2$. For $k=0,\ldots,\Delta/h$, let $a_k$ denote the sum of the state error and Jacobian error after $k$ Euler updates. The Lipschitz bounds on $f$ and $Df$ give
\[
a_{k+1}\le(1+Ch)a_k+Ch^2, \qquad a_0=0.
\]
Discrete Gronwall therefore yields
\[
\|\mathcal E_{h}-\phi_{\Delta}\|_{C^1} \le C\Delta h.
\]
Applying this estimate with $h=\lambda$ and $h=\tau$ gives
\[
\|\mathcal E_{\lambda} -\mathcal E_{\tau}\|_{C^1} \le C\Delta(\lambda+\tau).
\]

The decomposition
\[
\mathcal B_{j} = \bigl(\E_\xi\cF_{j} -\mathcal E_{\lambda}\bigr) + \bigl(\mathcal E_{\lambda} -\mathcal E_{\tau}\bigr)
\]
now gives
\[
\|\mathcal B_{j}\|_{C^1} \le C\Delta \bigl(\lambda+\tau+s_n\sqrt{\Delta}\bigr).
\]
Its derivative and value bounds prove (N3) and (N4).

Finally,
\[
s_n+d_n \le (1+\sqrt T)(s_n+\lambda+\tau).
\]
The path estimate follows from Theorem~\ref{thm:nonlinear-uniform}, after enlarging the constant $C_T$.
\end{proof}

The mean of a nonlinear fine propagator need not equal noise-free Euler propagation, despite the conditional unbiasedness of each learner update. This accounts for the term $s_n\sqrt{\Delta}$ in the mismatch bound. For fixed $T$, this term is absorbed into the learner-noise scale in the resulting error rate.

The bounded learner draws also yield conditional moment bounds of every fixed order.

\begin{lemma}\label{lem:Lr}
Under Assumption~\ref{ass:primitive}, fix $q\ge2$. The following statements hold with constants independent of $n$, $M$, $L$, $J$ and $B$.
\begin{enumerate}[label=(\alph*)]
\item For deterministic $x,y\in\R^p$,
\[
\begin{aligned}
\|\mathcal Z_{j}(x)-\mathcal Z_{j}(y)\|_{L^q_\xi} &\le C_q\sqrt{\Delta}\,s_n\|x-y\|,\\
\|\mathcal Z_{j}(x)\|_{L^q_\xi} &\le C_q\sqrt{\Delta}\,s_n.
\end{aligned}
\]

\item The conclusion of Lemma~\ref{lem:nonlin-predictable} holds in conditional $L^q$ when its conditional second-moment bound is replaced by
\[
\E(\|z_j\|^q\mid\mathcal H_{j}) \le C_q\Delta^{q/2}s_n^qa_j^q,
\]
and $\|\max_j a_j\|_{L^q_\xi}<\infty$.

\item There exists $C_{T,q}\ge1$, independent also of $K$, such that
\[
\left( \E_\xi \|\Theta^{(K,M)}-\Theta\|_\infty^q \right)^{1/q} \le C_{T,q}^{K+2}\eta_n^{K+1}, \qquad K\ge0.
\]
\end{enumerate}
The constants may depend on $q$, $T$, $p$, $m$, $C_f$, $L_\xi$ and $\|\beta_0\|$.
\end{lemma}

\begin{proof}
The centered, conditionally independent learner draws satisfy, by Burkholder's inequality,
\[
\begin{aligned}
\E_\xi\|\bar\xi_{\ell}\|^q &\le C_{q,m}B^{-q} \E_\xi \left( \sum_{r=1}^{B}\|\xi_{\ell,r}\|^2 \right)^{q/2}\\
&\le C_{q,m}L_\xi^qB^{-q/2}.
\end{aligned}
\]
Thus the batch-average moment constants are uniform in $B$.

Use the notation of the preceding proof. The martingale part of $E_\ell$ has increments $\lambda\sigma(X_\ell)\zeta_{\ell+1}$. Burkholder's inequality and Minkowski's inequality bound its maximal $L^q$ norm by
\[
C_q\lambda\sqrt{L/B} = C_qs_n\sqrt{\Delta}.
\]
The pathwise discrete Gronwall argument gives the same bound for $\max_{\ell\le L}\|E_\ell\|$. The Jacobian difference has a bounded noise coefficient and drift bounded by $C(\|\Gamma_\ell\|_F+\|E_\ell\|)$. The same argument yields
\[
\left\| \max_{\ell\le L}\|\Gamma_\ell\|_F \right\|_{L^q_\xi} \le C_qs_n\sqrt{\Delta}.
\]
Centering and the integral representation of the coupled fine-map difference prove part (a).

For part (b), retain the pathwise Gronwall estimate of Lemma~\ref{lem:nonlin-predictable}. Set $S=\max_{i\le M}\|\sum_{j<i}z_j\|$. Burkholder's inequality and Minkowski's inequality give
\[
\|S\|_{L^q_\xi}^2 \le C_q \left\|\sum_{j<M}\|z_j\|^2\right\|_{L^{q/2}_\xi} \le C_q\sum_{j<M}\|z_j\|_{L^q_\xi}^2 \le C_qTs_n^2 \left\|\max_{j<M}a_j\right\|_{L^q_\xi}^2.
\]
This proves the conditional $L^q$ forcing estimate.

For part (c), condition on $\mathcal H_{j}$ and apply part (a) to the centered difference in \eqref{eq:nonlin-exact}. Part (b) gives the correction recursion in $L^q$. The sequential fine path is bounded pathwise, so the predictor has the corresponding $L^q$ estimate. Iteration and the reconstruction inequality prove the result.
\end{proof}

\subsection{Randomized coordinate boosting}
\label{subsec:coordinate}

We now verify the propagator assumptions for the coordinate learner defined in Section~\ref{sec:setting}. Its affine structure permits a direct calculation of the centered fluctuation and the deterministic mismatch. Define $\pi_{\min}=\min_{1\le i\le p}\pi_i$.

\begin{assumption}\label{ass:bounded}
The coordinate probabilities are fixed, with $\pi_{\min}>0$, and
\[
\|H_n\|+\|b_n\|=O_\PP(1).
\]
\end{assumption}

For $R>0$, define the data-localization event
\[
\Omega_{n,R} = \{\|H_n\|+\|b_n\|\le R\}.
\]
For $x\in\R^p$, let $\bar x=(x^\top,1)^\top\in\R^{p+1}$. Define the augmented matrices $\mathsf L_n,\widehat{\mathsf L}_{\ell} \in\R^{(p+1)\times(p+1)}$ by
\[
\mathsf L_n =
\begin{pmatrix}
-H_n&b_n\\
0&0
\end{pmatrix},
\qquad \widehat{\mathsf L}_{\ell} =
\begin{pmatrix}
-\widehat D_{\ell}H_n&\widehat D_{\ell}b_n\\
0&0
\end{pmatrix}.
\]
Let $I_{p+1}$ be the identity matrix in the augmented space. For each subinterval $j$, define
\[
\mathsf P_{j,0} =I_{p+1}, \quad \mathsf P_{j,\ell+1} = \bigl( I_{p+1} +\lambda\widehat{\mathsf L}_{jL+\ell+1} \bigr)\mathsf P_{j,\ell}, \qquad 0\le\ell<L.
\]
Set
\[
\mathsf P_{j}=\mathsf P_{j,L}, \qquad \mathsf G=(I_{p+1}+\tau\mathsf L_n)^{J}, \qquad \mathsf Q_{j} = \mathsf P_{j}-\E_\xi\mathsf P_{j}.
\]
The matrices $\mathsf P_{j}$ and $\mathsf G$ represent the fine and coarse affine propagators.

We first estimate their centered fluctuation and deterministic difference.

\begin{lemma}\label{lem:propagator}
Under Assumption~\ref{ass:bounded}, on $\Omega_{n,R}$,
\[
\begin{aligned}
\E_\xi\|\mathsf Q_{j}\|_F^2 &\le C_R\Delta s_n^2,\\
\|\E_\xi\mathsf P_{j}-\mathsf G\|_F &\le C_R\Delta(\lambda+\tau).
\end{aligned}
\]
All fine partial products and coarse products over total time at most $T$ have norm bounded by $C_R$. The constant depends on $R$, $T$, $p$ and $\pi_{\min}$, but not on the discretization levels or $B$.
\end{lemma}

\begin{proof}
On $\Omega_{n,R}$,
\[
\|\widehat{\mathsf L}_{\ell}\| +\|\mathsf L_n\|\le C_R
\]
uniformly in the coordinate draws. The product bounds follow from $\|I_{p+1}+hA\|\le1+h\|A\|$ and the fact that the total elapsed time is at most $T$.

Set $\overline{\mathsf W}=I_{p+1}+\lambda\mathsf L_n$. Conditional independence and $\E_\xi\widehat D_{\ell}=I_p$ give
\[
\E_\xi\mathsf P_{j} = \overline{\mathsf W}^{L}.
\]
The exact telescoping identity is
\[
\mathsf Q_{j} = \lambda\sum_{\ell=1}^{L} \overline{\mathsf W}^{L-\ell} \bigl( \widehat{\mathsf L}_{jL+\ell}-\mathsf L_n \bigr) \mathsf P_{j,\ell-1}.
\]
Each summand is centered conditionally on the preceding local draws. Furthermore,
\[
\E_\xi \|\widehat{\mathsf L}_{\ell}-\mathsf L_n\|_F^2 \le \frac{C_R}{B},
\]
because $\widehat D_{\ell}$ averages $B$ independent coordinate matrices. Orthogonality and the uniform product bounds therefore yield
\[
\E_\xi\|\mathsf Q_{j}\|_F^2 \le C_RL\lambda^2/B = C_R\Delta s_n^2.
\]

For the deterministic difference, the one-step matrix-exponential error is bounded by $C_Rh^2$. Telescoping this error over a subinterval gives
\[
\begin{aligned}
\| \overline{\mathsf W}^{L} -e^{\Delta\mathsf L_n} \|_F &\le C_R\Delta\lambda,\\
\| \mathsf G-e^{\Delta\mathsf L_n} \|_F &\le C_R\Delta\tau.
\end{aligned}
\]
Their sum proves the bias estimate.
\end{proof}

These matrix estimates verify the hypotheses of the general error theorem.

\begin{theorem}\label{thm:parareal}
Under Assumption~\ref{ass:bounded}, there exists $C_R\ge1$ such that, on $\Omega_{n,R}$ and for every $K\ge0$,
\begin{equation}\label{eq:parareal-endpoint}
\left( \E_\xi\max_{0\le j\le M} \|U_{j}^{(K)}-\Theta(T_j)\|^2 \right)^{1/2} \le (C_R\eta_n)^{K+1},
\end{equation}
and
\begin{equation}\label{eq:parareal-uniform}
\left( \E_\xi \|\Theta^{(K,M)}-\Theta\|_\infty^2 \right)^{1/2} \le C_R(C_R\eta_n)^{K+1}.
\end{equation}
The constant may depend on $\|\beta_0\|$ in addition to the quantities in Lemma~\ref{lem:propagator}, but is independent of $n$, $M$, $L$, $J$, $B$ and $K$. Consequently, for each fixed $K$,
\[
\|\Theta^{(K,M)}-\Theta\|_\infty = O_\PP(\eta_n^{K+1}).
\]
\end{theorem}

\begin{proof}
The centered fine map is represented by the upper components of $\mathsf Q_{j}\bar x$. Since $\|\bar x\|\le1+\|x\|$, Lemma~\ref{lem:propagator} proves (N2). Applying the same matrix to $(x-y,0)^\top$ proves (N1). The deterministic matrix mismatch gives (N3) and (N4) with $d_n=\lambda+\tau$.

Let $G$ denote the linear part of $\cG$. The identity
\[
\mathsf G-I_{p+1} = \tau\sum_{k=0}^{J-1} \mathsf L_n(I_{p+1}+\tau\mathsf L_n)^k
\]
and the product bounds imply $\|G-I_p\|\le C_R\Delta$. The partial fine maps are uniformly Lipschitz. The matrix-product bounds also imply conditional square integrability at every deterministic starting state. Applied to the augmented initial vector $(\beta_0^\top,1)^\top$, they bound the complete sequential fine path by $C_R$. Thus all requirements of Assumption~\ref{ass:nonlin-propagator} hold on $\Omega_{n,R}$.

Theorem~\ref{thm:nonlinear-uniform} gives \eqref{eq:parareal-endpoint} and \eqref{eq:parareal-uniform}, after enlarging $C_R$. For the unconditional assertion, choose $R$ so that $\PP(\Omega_{n,R}^c)$ is small and apply conditional Markov's inequality on $\Omega_{n,R}$. This proves the stated $O_\PP$ bound.
\end{proof}

Both learner models give a common consequence after one correction.

\begin{corollary}\label{cor:one}
Under either Proposition~\ref{prop:primitive} or Theorem~\ref{thm:parareal},
\[
\|\Theta^{(1,M)}-\Theta\|_\infty = O_\PP\bigl((s_n+\lambda+\tau)^2\bigr).
\]
If, in addition, $\lambda/\tau\to0$ and $\tau=O(s_n)$, then
\[
\|\Theta^{(1,M)}-\Theta\|_\infty = O_\PP(s_n^2) = o_\PP(s_n).
\]
\end{corollary}

\begin{proof}
The first bound is the case $K=1$ of the two error estimates. Under the additional conditions,
\[
\frac{\lambda}{s_n} = \frac{\lambda}{\tau}\frac{\tau}{s_n} \longrightarrow0,
\]
so $\eta_n=O(s_n)$. Since $s_n^2=\lambda/B\le\lambda=o(s_n)$, also $s_n\to0$, which gives the last equality.
\end{proof}

\section{Asymptotic distribution of the Parareal error}
\label{sec:second}

We identify the limiting error for every fixed correction level $K<M$. The number of subintervals remains fixed. For $K\ge M$, Proposition~\ref{prop:finite-exactness} gives exact agreement with the sequential trajectory.

\subsection{Limit theorem}

We retain the fine maps, the Euler coarse map and the continuous reconstruction defined in Section~\ref{sec:setting}. In particular,
\[
\Delta=\frac TM,\qquad T_j=j\Delta,\qquad s_n=\sqrt{\frac{\lambda}{B}}.
\]
We write $C_b^2$ for the space of twice continuously differentiable fields with bounded values and bounded first and second derivatives. For a compact set $D\subset\R^p$, the $C^k(D)$ norm is the sum of the supremum norms of the derivatives through order $k$. We omit $D$ when the supremum is taken over $\R^p$.

\begin{assumption}\label{ass:multi}
The learner has the representation
\[
\widehat F(x,\xi)=f_n(x)+\sigma_n(x)\xi,
\]
where $f_n:\R^p\to\R^p$ and $\sigma_n:\R^p\to\R^{p\times m}$ are $\mathcal D_n$-measurable $C_b^2$ fields and $m\ge1$ is fixed. There are a constant $C$ and events $\Omega_n\in\mathcal D_n$ with $\PP(\Omega_n)\to1$ such that
\[
\|f_n\|_{C^2}+\|\sigma_n\|_{C^2}\le C \quad\text{on }\Omega_n.
\]

For each $n$, the array $(\xi_{\ell,r})_{1\le\ell\le N,\;1\le r\le B}$ is independent of $\mathcal D_n$ and has independent entries with a common law that does not depend on $n$. A vector $\xi$ with this law satisfies
\[
\E\xi=0,\qquad \E(\xi\xi^\top)=I_m,\qquad \|\xi\|\le L_\xi\quad\text{almost surely}.
\]

There are deterministic $C_b^2$ fields $f$ and $\sigma$ such that, for every compact $D\subset\R^p$,
\[
\|f_n-f\|_{C^2(D)} +\|\sigma_n-\sigma\|_{C^2(D)} \xrightarrow{\PP}0.
\]
Moreover,
\[
\frac{\lambda}{\tau}\longrightarrow0, \qquad \frac{\tau}{s_n}\longrightarrow c \quad\text{for some }c\in[0,\infty].
\]
\end{assumption}

Let $\phi_{n,t}$ and $\phi_t$ be the flows generated by $f_n$ and $f$, respectively, and define
\[
\theta_n(t)=\phi_{n,t}(\beta_0),\qquad \theta(t)=\phi_t(\beta_0),\qquad A(t)=Df(\theta(t)).
\]
These flows are defined for all finite times under Assumption~\ref{ass:multi}. For $0\le u\le t\le T$, let $\Phi(t,u)$ solve
\[
\partial_t\Phi(t,u)=A(t)\Phi(t,u),\qquad \Phi(u,u)=I_p,
\]
and put
\[
\Phi_j=\Phi(T_{j+1},T_j),\qquad 0\le j<M.
\]

For $1\le\ell\le N$, define
\[
\bar\xi_{\ell} =\frac1{B}\sum_{r=1}^{B}\xi_{\ell,r}.
\]
Let $W_n$ be the piecewise-linear path with values
\[
W_n(\ell\lambda) =\frac{\lambda}{s_n} \sum_{k=1}^{\ell}\bar\xi_{k}, \qquad 0\le\ell\le N.
\]
Each increment over a fine step has covariance $\lambda I_m$.

Let $W$ be an $m$-dimensional Brownian motion. The first-order learner fluctuation is the Gaussian process $V$ satisfying
\begin{equation}\label{eq:multi-V}
dV(t)=A(t)V(t)\,dt+\sigma(\theta(t))\,dW(t), \qquad V(0)=0.
\end{equation}

To describe the Euler bias, define
\[
h(x)=Df(x)f(x).
\]
Its derivative is
\[
Dh(x)[v] =D^2f(x)[v,f(x)]+Df(x)Df(x)v, \qquad x,v\in\R^p.
\]
Let $b_{\mathrm E}$ solve
\begin{equation}\label{eq:critical-b}
\dot b_{\mathrm E}(t) =A(t)b_{\mathrm E}(t)-\frac12h(\theta(t)), \qquad b_{\mathrm E}(0)=0.
\end{equation}
The corresponding bias over one subinterval is
\begin{equation}\label{eq:propagator-euler-bias}
\mathcal E_\Delta(x) =-\frac12\int_0^\Delta D\phi_{\Delta-u}(\phi_u(x))h(\phi_u(x))\,du.
\end{equation}
Thus
\[
b_{\mathrm E}(T_{j+1}) =\Phi_jb_{\mathrm E}(T_j) +\mathcal E_\Delta(\theta(T_j)).
\]
We write
\[
\mathsf E_j=D\mathcal E_\Delta(\theta(T_j)).
\]

The random sensitivity of a fine propagator is obtained by differentiating its first-order fluctuation. Define
\[
\mathcal K_\Delta(u,x) =D\phi_{\Delta-u}(\phi_u(x))\sigma(\phi_u(x)), \qquad 0\le u\le\Delta.
\]
For $0\le j<M$, set
\begin{align}
\eta_j &=\int_0^\Delta \mathcal K_\Delta(u,\theta(T_j))\,dW(T_j+u), \label{eq:propagator-gaussian-increment}\\
\mathcal R_j[v] &=\int_0^\Delta D_x\mathcal K_\Delta(u,\theta(T_j))[v]\, dW(T_j+u), \qquad v\in\R^p. \label{eq:propagator-random-sensitivity}
\end{align}
Here $\mathcal R_j$ is a random linear map from $\R^p$ to $\R^p$. The pairs $(\eta_j,\mathcal R_j)$ are centered Gaussian and independent across subintervals. Moreover,
\[
V(T_{j+1})=\Phi_jV(T_j)+\eta_j.
\]

For endpoint values $z_0,\ldots,z_M\in\R^p$, define the continuous path $\mathcal I_M(z)$ by
\begin{equation}\label{eq:limit-reconstruction}
\mathcal I_M(z)(T_j+u) = \Phi(T_j+u,T_j)z_j +\frac u\Delta\{z_{j+1}-\Phi_jz_j\}, \qquad 0\le u\le\Delta.
\end{equation}
It takes the value $z_j$ at $T_j$.

\begin{theorem}\label{thm:parareal-limit}
Under Assumption~\ref{ass:multi}, suppose first that $c<\infty$. Define the level-zero endpoint values by
\[
z_j^{0,c}=c\,b_{\mathrm E}(T_j)-V(T_j), \qquad 0\le j\le M,
\]
and the subsequent levels by
\begin{equation}\label{eq:general-limit-recursion}
z_0^{K+1,c}=0, \quad z_{j+1}^{K+1,c} =\Phi_jz_j^{K+1,c} +(\mathcal R_j-c\mathsf E_j)[z_j^{K,c}], \quad 0\le j<M,\quad K\ge0.
\end{equation}
Let
\[
Z^{K,c}=\mathcal I_M(z^{K,c}).
\]
Then, jointly in continuous-path spaces,
\begin{equation}\label{eq:multivariate-joint}
\left( W_n,\, \frac{\Theta-\theta_n}{s_n},\, \left( \frac{\Theta^{(K,M)}-\Theta}{s_n^{K+1}} \right)_{K=0}^{M-1} \right) \Rightarrow \left( W,V,(Z^{K,c})_{K=0}^{M-1} \right).
\end{equation}
The conditional law of the vector on the left, given $\mathcal D_n$, also converges weakly in probability to the same nonrandom limit law.

If $c=\infty$, define
\[
\bar z_j^0=b_{\mathrm E}(T_j),\qquad 0\le j\le M,
\]
and
\begin{equation}\label{eq:general-deterministic-recursion}
\begin{aligned}
\bar z_0^{K+1}&=0,\\
\bar z_{j+1}^{K+1} &=\Phi_j\bar z_j^{K+1}-\mathsf E_j[\bar z_j^K], \quad 0\le j<M,\quad K\ge0.
\end{aligned}
\end{equation}
For every $0\le K<M$,
\begin{equation}\label{eq:coarse-dominant-limit}
\frac{\Theta^{(K,M)}-\Theta}{\tau^{K+1}} \xrightarrow{\PP} \mathcal I_M(\bar z^K) \quad\text{in }C([0,T];\R^p).
\end{equation}
For $K\ge M$, the path error is identically zero.
\end{theorem}

We first establish the propagator expansions used in the proof.

\begin{lemma}\label{lem:propagator-expansions}
Under Assumption~\ref{ass:multi}, write
\[
P_{j}=U_{j}^{(0)},\qquad \theta_{n,j}=\theta_n(T_j),\qquad \delta_{j}=\mathcal F_{j}-\mathcal G.
\]
For every compact $D\subset\R^p$,
\begin{equation}\label{eq:coarse-C1-expansion}
\left\| \frac{\mathcal G-\phi_{n,\Delta}}{\tau} -\mathcal E_\Delta \right\|_{C^1(D)} \xrightarrow{\PP}0,
\end{equation}
and
\begin{equation}\label{eq:predictor-expansion}
\max_{j\le M} \left\| \frac{P_{j}-\theta_{n,j}}{\tau} -b_{\mathrm E}(T_j) \right\| \xrightarrow{\PP}0.
\end{equation}
Furthermore,
\begin{equation}\label{eq:propagator-C2-consistency}
\begin{aligned}
\max_{j<M}\sup_{0\le u\le\Delta} \|\mathcal F_{j,u}-\phi_{n,u}\|_{C^2(D)} &\xrightarrow{\PP}0,\\
\|\mathcal G-\phi_{n,\Delta}\|_{C^2(D)} &\xrightarrow{\PP}0.
\end{aligned}
\end{equation}
In particular,
\begin{equation}\label{eq:second-derivative-consistency}
\max_{j<M}\|D^2\delta_{j}\|_{C^0(D)} \xrightarrow{\PP}0.
\end{equation}

If $c<\infty$, then jointly with $W_n$,
\begin{equation}\label{eq:fine-propagator-derivative-limit}
\left( \frac{\Theta-\theta_n}{s_n},\, \left( \frac{ D\mathcal F_{j}(\theta_{n,j}) -D\phi_{n,\Delta}(\theta_{n,j}) }{s_n} \right)_{j=0}^{M-1} \right) \Rightarrow \left(V,(\mathcal R_j)_{j=0}^{M-1}\right).
\end{equation}
This joint convergence also holds conditionally on $\mathcal D_n$, weakly in probability.

For all values of $c$,
\begin{equation}\label{eq:propagator-first-order-bounds}
\begin{aligned}
\|\Theta-\theta_n\|_\infty &=O_\PP(s_n+\lambda),\\
\max_{j<M} \left\| D\mathcal F_{j}(\theta_{n,j}) -D\phi_{n,\Delta}(\theta_{n,j}) \right\| &=O_\PP(s_n+\lambda).
\end{aligned}
\end{equation}
\end{lemma}

\begin{proof}
For each fixed training sample, the coefficient norms are finite. Since $\|\bar\xi_{\ell}\|\le L_\xi$,
\[
\|\mathcal F_{j}(x)-x\| \le \Delta\{\|f_n\|_\infty +L_\xi\|\sigma_n\|_\infty\}.
\]
Thus every fine propagator has finite conditional moments at every deterministic starting state.

On $\Omega_n$, the preceding bound and its coarse counterpart give
\[
\|\mathcal F_{j}(x)-x\|\le C\Delta,\qquad \|\mathcal G(x)-x\|\le C\Delta.
\]
Consequently,
\[
\|U_{j+1}^{(K+1)}-U_{j}^{(K+1)}\| \le C\Delta.
\]
Indeed, this increment is the sum of $\mathcal G(U_{j}^{(K+1)})-U_{j}^{(K+1)}$ and $\delta_{j}(U_{j}^{(K)})$. All Parareal endpoints are therefore bounded by $\|\beta_0\|+CT$ on $\Omega_n$. Their intermediate fine and coarse trajectories are bounded as well.

To prove the convergence assertions, the subsequence principle allows us to work with deterministic coefficient sequences converging in $C^2$ on compact sets and satisfying the uniform bounds. Such sequences arise by conditioning along almost surely convergent training-data subsequences on which $\Omega_n$ eventually occurs.

Fix a compact set $D$ of starting states. Choose a closed ball $D_*$ containing the intermediate fine, coarse and mean-flow trajectories from $D$. The coefficient bounds make this choice independent of $n$. Convergence in $C^2(D_*)$ provides a common modulus of continuity $\omega$, with $\omega(r)\to0$ as $r\downarrow0$, for the second derivatives of $f_n$ and $\sigma_n$ on $D_*$.

We first derive the Euler bias together with its derivative. For a matrix $\mathsf J\in\R^{p\times p}$, define the augmented vector field
\[
\mathfrak f_n(x,\mathsf J) =\bigl(f_n(x),Df_n(x)\mathsf J\bigr).
\]
This field is $C^1$. Its exact trajectory started at $(x,I_p)$ is
\[
\bigl(\phi_{n,t}(x),D\phi_{n,t}(x)\bigr),
\]
and its Euler scheme is precisely the coarse recursion together with its Jacobian recursion.

On the bounded sets reached by these trajectories, the local exact-flow defect of the augmented Euler scheme is
\[
\frac{h^2}{2} D\mathfrak f_n(y)\mathfrak f_n(y)+o(h^2),
\]
uniformly in $n$ and $y$. The uniform remainder follows from the common modulus of continuity of $D\mathfrak f_n$. Propagating the local defects and linearizing the accumulated error shows that the global error divided by $h$ converges to the linear variational equation with forcing
\[
-\frac12D\mathfrak f_n(y(t))\mathfrak f_n(y(t)).
\]
Its state component is $\mathcal E_{n,t}(x)$, where $\mathcal E_{n,t}$ is defined by \eqref{eq:propagator-euler-bias} with the fields indexed by $n$ and $\Delta$ replaced by $t$. Its Jacobian component is $D\mathcal E_{n,t}(x)$, as follows by differentiating the state-bias equation.

Taking $h=\tau$ and $t=\Delta$ proves
\[
\left\| \frac{\mathcal G-\phi_{n,\Delta}}{\tau} -\mathcal E_{n,\Delta} \right\|_{C^1(D)} \longrightarrow0.
\]
Continuous dependence of the flow and its variational equations gives $\mathcal E_{n,\Delta}\to\mathcal E_\Delta$ in $C^1(D)$. This proves \eqref{eq:coarse-C1-expansion}. Applying the state component over $[0,T_j]$ gives \eqref{eq:predictor-expansion}, since $\mathcal E_{T_j}(\beta_0)=b_{\mathrm E}(T_j)$.

We next compare a fine propagator with its mean flow. Let $X_\ell^x$ be the fine trajectory started at $x$, and set
\[
\vartheta_\ell=\phi_{n,\ell\lambda}(x),\qquad \mathsf J_\ell=D_xX_\ell^x,\qquad \mathsf J^\phi_\ell=D\phi_{n,\ell\lambda}(x).
\]
Define
\[
E_\ell=X_\ell^x-\vartheta_\ell,\qquad \Gamma_\ell=\mathsf J_\ell-\mathsf J^\phi_\ell.
\]
The state and Jacobian recursions have martingale increments with accumulated conditional second moments bounded by $C\Delta s_n^2$. Their drift differences are bounded by $C(\|E_\ell\|+\|\Gamma_\ell\|)$, and the deterministic local flow defects are bounded by $C\lambda^2$. Doob's inequality and discrete Gronwall therefore give
\begin{equation}\label{eq:local-state-jacobian-estimate}
\left( \E_\xi\max_{\ell\le L} \{\|E_\ell\|+\|\Gamma_\ell\|\}^2 \right)^{1/2} \le C(s_n+\lambda),
\end{equation}
uniformly for $x\in D$. The same argument over the full fine grid proves \eqref{eq:propagator-first-order-bounds}.

For the second derivatives, write
\[
\mathsf Q_\ell[v,w]=D_x^2X_\ell^x[v,w], \qquad \mathsf Q^\phi_\ell[v,w] =D^2\phi_{n,\ell\lambda}(x)[v,w], \qquad v,w\in\R^p.
\]
Both initial values are zero. With $\zeta_{\ell+1}=\bar\xi_{jL+\ell+1}$,
\[
\begin{aligned}
\mathsf Q_{\ell+1}[v,w] ={}&\mathsf Q_\ell[v,w] +\lambda Df_n(X_\ell^x)\mathsf Q_\ell[v,w]\\
&+\lambda D^2f_n(X_\ell^x) [\mathsf J_\ell v,\mathsf J_\ell w]\\
&+\lambda\Bigl\{ D\sigma_n(X_\ell^x)[\mathsf Q_\ell[v,w]] +D^2\sigma_n(X_\ell^x) [\mathsf J_\ell v,\mathsf J_\ell w] \Bigr\}\zeta_{\ell+1}.
\end{aligned}
\]
The bounded derivatives and discrete Gronwall give uniform pathwise bounds for $\mathsf J_\ell$ and $\mathsf Q_\ell$. Denote the accumulated martingale in this recursion by $M_\ell^Q$. Componentwise application of Doob's inequality yields
\[
\E_\xi\max_{\ell\le L}\|M_\ell^Q\|^2 \le C\Delta s_n^2.
\]

Comparison with the second variational equation of the mean flow gives
\[
\begin{aligned}
\max_{\ell\le L} \|\mathsf Q_\ell-\mathsf Q^\phi_\ell\| \le C\Bigl\{& \max_{\ell\le L}\|M_\ell^Q\| +\max_{\ell\le L}\|E_\ell\| +\max_{\ell\le L}\|\Gamma_\ell\|\\
&+\omega\!\left(\max_{\ell\le L}\|E_\ell\|\right) +\lambda+\omega(C\lambda)\Bigr\}.
\end{aligned}
\]
The last two terms bound the accumulated deterministic consistency error of the second variational equation. Every term on the right tends to zero in probability.

Uniformity in the starting state follows from a pathwise equicontinuity estimate. For the same stored innovations and two starting states $x,y\in D$, the differentiated recursions give
\[
\max_{\ell\le L} \|D_x^2X_\ell^x-D_x^2X_\ell^y\| \le C\{\|x-y\|+\omega(C\|x-y\|)\}.
\]
A finite-net argument therefore makes the convergence uniform on $D$. The coarse second derivatives satisfy the same comparison with step $\tau$ and without the martingale term. Linear interpolation of the fine maps preserves the convergence, because the mean-flow derivatives are uniformly continuous in time. This proves \eqref{eq:propagator-C2-consistency}, and subtraction gives \eqref{eq:second-derivative-consistency}.

It remains to identify the normalized first-order fluctuations when $c<\infty$. The mesh condition gives
\[
\frac{\lambda}{s_n} =\frac{\lambda}{\tau}\frac{\tau}{s_n} \longrightarrow0.
\]
The increments of $W_n$ are centered, have covariance $\lambda I_m$, and are bounded by $L_\xi\lambda/s_n$. Their predictable quadratic variation converges to $tI_m$, so the martingale functional central limit theorem gives $W_n\Rightarrow W$.

By \eqref{eq:local-state-jacobian-estimate}, $E_\ell$ and $\Gamma_\ell$ are $O_\PP(s_n)$ uniformly over the local grid. The state drift remainder, after division by $s_n$, is bounded by
\[
C s_n^{-1}\max_{\ell\le L}\|E_\ell\|^2 =O_\PP(s_n).
\]
For the Jacobian drift, Taylor expansion gives
\[
\begin{aligned}
&\bigl\| Df_n(X_\ell^x)\mathsf J_\ell -Df_n(\vartheta_\ell)\mathsf J^\phi_\ell -Df_n(\vartheta_\ell)\Gamma_\ell -D^2f_n(\vartheta_\ell)[E_\ell,\mathsf J^\phi_\ell] \bigr\| \\
&\qquad\le C\|E_\ell\|\|\Gamma_\ell\| +C\omega(\|E_\ell\|)\|E_\ell\|.
\end{aligned}
\]
Here the second derivative acts on each column of $\mathsf J^\phi_\ell$. After summation with weights $\lambda/s_n$, this remainder tends to zero in probability.

The changes in the state and Jacobian noise coefficients are bounded by $C(\|E_\ell\|+\|\Gamma_\ell\|)$. After normalization, their accumulated quadratic variation is bounded by
\[
C\lambda\sum_{\ell<L} (\|E_\ell\|+\|\Gamma_\ell\|)^2.
\]
Its expectation tends to zero by \eqref{eq:local-state-jacobian-estimate}. The normalized deterministic flow defects also vanish, since their accumulated size is $O(\lambda/s_n)$.

Applying the joint martingale limit argument to the global state and the local state and Jacobian recursions therefore identifies their limits. On subinterval $j$, the local limits satisfy
\begin{align*}
d\mathcal V_j(t) &=A(t)\mathcal V_j(t)\,dt +\sigma(\theta(t))\,dW(t),\\
d\mathcal J_j(t)[v] &=A(t)\mathcal J_j(t)[v]\,dt + D^2f(\theta(t))[\mathcal V_j(t),\Phi(t,T_j)v]\,dt\\
&\quad + D\sigma(\theta(t))[\Phi(t,T_j)v]\,dW(t),
\end{align*}
with zero initial values at $T_j$. The global state limit is $V$ in \eqref{eq:multi-V}.

Variation of constants gives the kernel $\mathcal K_\Delta$ for $\mathcal V_j(T_{j+1})$. Its bounded deterministic derivative permits differentiation in $L^2$, yielding
\[
\mathcal J_j(T_{j+1})[v]=\mathcal R_j[v].
\]
All recursions use the same innovation array, so the joint limit uses the same Brownian motion. This proves \eqref{eq:fine-propagator-derivative-limit}.

The preceding argument applies conditionally along every almost surely convergent coefficient subsequence. It therefore gives conditional weak convergence in probability. Taking expectations gives the unconditional convergence.
\end{proof}

\begin{proof}[Proof of Theorem~\ref{thm:parareal-limit}]
Set
\[
X_{n,j}=\Theta(T_j),\qquad e_{j}^K=U_{j}^{(K)}-X_{n,j},\qquad \Phi_{n,j}=D\phi_{n,\Delta}(\theta_{n,j}).
\]
For every $K\ge0$, the Parareal recursion gives exactly
\begin{equation}\label{eq:general-error-recursion}
e_{j+1}^{K+1} = \mathcal G(X_{n,j}+e_{j}^{K+1}) -\mathcal G(X_{n,j}) +\delta_{j}(X_{n,j}+e_{j}^K) -\delta_{j}(X_{n,j}).
\end{equation}

Suppose first that $c<\infty$. Since $\lambda=o(s_n)$ and $\tau=O(s_n)$, the bounds of Section~\ref{sec:bounds} give
\[
\max_{j\le M}\|e_{j}^K\| =O_\PP(s_n^{K+1})
\]
for every fixed $K$. Define
\[
\varepsilon_{n,j}^K=s_n^{-(K+1)}e_{j}^K.
\]
At level zero, Lemma~\ref{lem:propagator-expansions} gives
\[
(\varepsilon_{n,j}^0)_{j=0}^{M} \Rightarrow (c\,b_{\mathrm E}(T_j)-V(T_j))_{j=0}^{M},
\]
jointly with $W_n$, the fine fluctuation and the normalized propagator derivatives.

Write the two terms in \eqref{eq:general-error-recursion} as
\[
\begin{aligned}
\mathcal G(X_{n,j}+e_{j}^{K+1}) -\mathcal G(X_{n,j}) &=\Phi_{n,j}e_{j}^{K+1}+r_{n,j}^{G,K},\\
\delta_{j}(X_{n,j}+e_{j}^K) -\delta_{j}(X_{n,j}) &=D\delta_{j}(\theta_{n,j})e_{j}^K +r_{n,j}^{\delta,K}.
\end{aligned}
\]
The pathwise bounds in Lemma~\ref{lem:propagator-expansions} place the relevant states in a fixed closed ball $D$ on $\Omega_n$. Uniform boundedness of the second coarse derivative and the $C^1$ Euler error give
\[
\|r_{n,j}^{G,K}\| \le C\{\tau+\|X_{n,j}-\theta_{n,j}\| +\|e_{j}^{K+1}\|\} \|e_{j}^{K+1}\| =O_\PP(s_n^{K+3}).
\]
For the defect term, the mean-value formula gives
\[
\begin{aligned}
\|r_{n,j}^{\delta,K}\| &\le C\|D^2\delta_{j}\|_{C^0(D)} \{\|X_{n,j}-\theta_{n,j}\|+\|e_{j}^K\|\} \|e_{j}^K\|\\
&=o_\PP(s_n^{K+2}).
\end{aligned}
\]
These estimates hold also for $K=0$.

Consequently,
\begin{equation}\label{eq:normalized-general-recursion}
\varepsilon_{n,j+1}^{K+1} =\Phi_{n,j}\varepsilon_{n,j}^{K+1} +\frac{D\delta_{j}(\theta_{n,j})}{s_n} [\varepsilon_{n,j}^K] +o_\PP(1).
\end{equation}
The propagator expansions imply
\[
\Phi_{n,j}\xrightarrow{\PP}\Phi_j, \qquad \left( \frac{D\delta_{j}(\theta_{n,j})}{s_n} \right)_{j=0}^{M-1} \Rightarrow (\mathcal R_j-c\mathsf E_j)_{j=0}^{M-1},
\]
jointly with the level-zero errors. Induction over the correction levels and the fixed number of subintervals proves joint convergence of the endpoint errors to \eqref{eq:general-limit-recursion}.

The continuous reconstruction satisfies
\begin{align*}
\Theta^{(K,M)}(T_j+u)-\Theta(T_j+u) = & \mathcal F_{j,u}(X_{n,j}+e_{j}^K) -\mathcal F_{j,u}(X_{n,j})\\
&\quad+\frac u\Delta \left\{ e_{j+1}^K -\mathcal F_{j}(X_{n,j}+e_{j}^K) +\mathcal F_{j}(X_{n,j}) \right\}.
\end{align*}
Using \eqref{eq:propagator-C2-consistency}, uniformly in $u$, we obtain
\[
\Theta^{(K,M)}(T_j+u)-\Theta(T_j+u) = D\phi_{n,u}(\theta_{n,j})e_{j}^K +\frac u\Delta \{e_{j+1}^K-\Phi_{n,j}e_{j}^K\} +o_\PP(s_n^{K+1}).
\]
After normalization, this is precisely the interpolation \eqref{eq:limit-reconstruction}. It proves the joint path convergence. The conditional convergence follows from the conditional propagator expansions and the same finite recursions.

Suppose now that $c=\infty$. Then $s_n=o(\tau)$ and $\lambda=o(\tau)$. The Section~\ref{sec:bounds} estimates give
\[
\max_{j\le M}\|e_{j}^K\| =O_\PP(\tau^{K+1}).
\]
Moreover,
\[
\frac{e_{j}^0}{\tau} \xrightarrow{\PP}b_{\mathrm E}(T_j), \qquad \frac{D\delta_{j}(\theta_{n,j})}{\tau} \xrightarrow{\PP}-\mathsf E_j.
\]
The same Taylor estimates, with $\tau$ in place of $s_n$, give \eqref{eq:general-deterministic-recursion}. The reconstruction argument proves \eqref{eq:coarse-dominant-limit}. Finally, Proposition~\ref{prop:finite-exactness} gives zero error for $K\ge M$.
\end{proof}

The conditional statement separates learner randomness from sampling fluctuations.

\begin{corollary}\label{cor:conditional}
Under Assumption~\ref{ass:multi}, with $c<\infty$, let $Y_n$ be a $\mathcal D_n$-measurable random vector such that $Y_n\Rightarrow Y$. Then
\[
\left( Y_n,\, \frac{\Theta-\theta_n}{s_n},\, \left( \frac{\Theta^{(K,M)}-\Theta}{s_n^{K+1}} \right)_{K=0}^{M-1} \right) \Rightarrow \left(Y,V,(Z^{K,c})_{K=0}^{M-1}\right),
\]
where $Y$ is independent of the limiting learner processes.
\end{corollary}

\begin{proof}
Condition on $\mathcal D_n$ in the expectation of a product of bounded Lipschitz test functions. By Theorem~\ref{thm:parareal-limit}, the conditional expectation involving the learner processes converges in probability to its nonrandom limiting expectation. Boundedness permits passage to expectations and gives the stated product law.
\end{proof}

\subsection{Distribution of the limit}

The recursion describes the limit through the finite Gaussian collection $(\eta_j,\mathcal R_j)_{j=0}^{M-1}$. It also identifies its Wiener-chaos components.

\begin{proposition}\label{prop:chaos-orders}
For every $0\le K<M$, there are processes $Z_r^K$, $0\le r\le K+1$, whose scalar coordinates belong to the Wiener chaos of order $r$ with respect to $W$, such that
\begin{equation}\label{eq:general-chaos-decomposition}
Z^{K,c} =\sum_{r=0}^{K+1}c^{K+1-r}Z_r^K.
\end{equation}
The component of order zero is deterministic and equals
\[
Z_0^K=\mathcal I_M(\bar z^K).
\]
In particular,
\[
Z^{K,0}=Z_{K+1}^K.
\]
For every $t\in[0,T]$,
\begin{align}
\E Z^{K,c}(t) &=c^{K+1}Z_0^K(t), \label{eq:general-limit-mean}\\
\operatorname{Cov}(Z^{K,c}(t)) &=\sum_{r=1}^{K+1} c^{2(K+1-r)}\operatorname{Cov}(Z_r^K(t)). \label{eq:general-limit-covariance}
\end{align}
\end{proposition}

\begin{proof}
At level zero, define endpoint components
\[
z_{j,0}^0=b_{\mathrm E}(T_j),\qquad z_{j,1}^0=-V(T_j).
\]
Set $z_{j,r}^K=0$ when $r<0$ or $r>K+1$. For subsequent levels, set $z_{0,r}^{K+1}=0$ and define
\[
z_{j+1,r}^{K+1} =\Phi_jz_{j,r}^{K+1} +\mathcal R_j[z_{j,r-1}^K] -\mathsf E_j[z_{j,r}^K].
\]
Every component at node $j$ depends only on Brownian increments before $T_j$. Multiplication by $\mathcal R_j$ therefore raises its chaos order by one: the current Gaussian propagator is independent of the preceding propagators, so the product has no contractions between those factors. The deterministic matrices preserve the chaos order.

Collecting powers of $c$ in \eqref{eq:general-limit-recursion} proves
\[
z_j^{K,c}=\sum_{r=0}^{K+1}c^{K+1-r}z_{j,r}^K.
\]
Applying the deterministic linear map $\mathcal I_M$ gives \eqref{eq:general-chaos-decomposition}, with $Z_r^K=\mathcal I_M(z_{\cdot,r}^K)$. The component of order zero satisfies \eqref{eq:general-deterministic-recursion}. Positive-order components are centered, and distinct Wiener chaoses are orthogonal. This proves the mean and covariance formulas.
\end{proof}

For one correction, the second-chaos component has a quadratic spectral representation.

\begin{corollary}\label{cor:spectral}
For every $\ell\in\R^p$, there are real numbers $\mu_1(\ell),\ldots,\mu_d(\ell)$, with $d\le Mp(p+1)$, such that
\begin{equation}\label{eq:multi-spectral}
\ell^\top Z^{1,0}(T) \stackrel{d}{=} \sum_{i=1}^{d}\mu_i(\ell)(G_i^2-1), \qquad G_i\ \text{i.i.d. }N(0,1).
\end{equation}
Moreover,
\[
\operatorname{Var}(\ell^\top Z^{1,0}(T)) =2\sum_{i=1}^{d}\mu_i(\ell)^2, \qquad \kappa_4(\ell^\top Z^{1,0}(T)) =48\sum_{i=1}^{d}\mu_i(\ell)^4,
\]
where $\kappa_4$ denotes the fourth cumulant. Every such contrast with positive variance is non-Gaussian.
\end{corollary}

\begin{proof}
The endpoint values $V(T_j)$ are deterministic linear combinations of $\eta_0,\ldots,\eta_{j-1}$. Hence the recursion expresses the contrast as a centered quadratic form in the Gaussian vector formed by the coordinates of $\eta_j$ and the matrix entries of $\mathcal R_j$, over all propagators. This vector has dimension $Mp(p+1)$.

Represent it in $d$ independent standard Gaussian coordinates, where $d$ is the rank of its covariance matrix. The contrast is then a quadratic form with a symmetric matrix of trace zero. Diagonalization gives \eqref{eq:multi-spectral}, with $\mu_i(\ell)$ the eigenvalues of that matrix. The variance and fourth-cumulant identities follow from those of $G_i^2-1$. Positive variance implies a strictly positive fourth cumulant.
\end{proof}

\subsection{Correction of the leading error}

The leading error at level $K\ge1$ can be estimated using two successive correction levels and derivatives of the propagators.

For $1\le K<M$ and $0\le j<M$, set
\[
v_{j}^K=U_{j}^{(K-1)}-U_{j}^{(K)}.
\]
Define $\widehat e_{0}^K=0$ and
\begin{equation}\label{eq:propagator-error-estimator}
\begin{aligned}
\widehat e_{j+1}^K ={}&D\mathcal G(U_{j}^{(K-1)})\widehat e_{j}^K\\
&+\left\{ D\mathcal F_{j}(U_{j}^{(K-1)}) -D\mathcal G(U_{j}^{(K-1)}) \right\}v_{j}^K.
\end{aligned}
\end{equation}
For $0\le u\le\Delta$, let
\begin{equation}\label{eq:propagator-path-estimator}
\begin{aligned}
\widehat e^K(T_j+u) ={}&D\mathcal F_{j,u}(U_{j}^{(K-1)})\widehat e_{j}^K\\
&+\frac u\Delta \left\{ \widehat e_{j+1}^K -D\mathcal F_{j}(U_{j}^{(K-1)})\widehat e_{j}^K \right\}.
\end{aligned}
\end{equation}
This path is continuous and satisfies $\widehat e^K(T_j)=\widehat e_{j}^K$. Define
\[
\widetilde\Theta^{(K,M)} =\Theta^{(K,M)}-\widehat e^K.
\]

\begin{proposition}\label{prop:error-correction}
Under Assumption~\ref{ass:multi}, with $c<\infty$, every fixed $1\le K<M$ satisfies
\[
\|\widetilde\Theta^{(K,M)}-\Theta\|_\infty =o_\PP(s_n^{K+1}).
\]
\end{proposition}
\begin{proof}
Retain the endpoint errors from the proof of Theorem~\ref{thm:parareal-limit}. We have
\[
e_{j}^{K-1}=v_{j}^K+e_{j}^K,\qquad e_{j}^{K-1}=O_\PP(s_n^K),\qquad e_{j}^K=O_\PP(s_n^{K+1}).
\]
The propagator expansions and \eqref{eq:second-derivative-consistency} also give
\[
D\delta_{j}(U_{j}^{(K-1)})=O_\PP(s_n).
\]

Taylor expansion of the exact error recursion around $U_{j}^{(K-1)}$ yields
\[
e_{j+1}^K = D\mathcal G(U_{j}^{(K-1)})e_{j}^K +D\delta_{j}(U_{j}^{(K-1)})v_{j}^K +o_\PP(s_n^{K+1}).
\]

The remainder $r_{n,j}^K$ in this expansion satisfies
\[
\begin{aligned}
\|r_{n,j}^K\| \le{}& C(\|e_{j}^{K-1}\|+\|e_{j}^K\|) \|e_{j}^K\| +\frac12\|D^2\delta_{j}\|_{C^0(D)} \|e_{j}^{K-1}\|^2\\
&+\|D\delta_{j}(U_{j}^{(K-1)})\| \|e_{j}^K\|,
\end{aligned}
\]
on a fixed closed ball $D$ containing the relevant states on $\Omega_n$. The three terms are respectively $O_\PP(s_n^{2K+1})$, $o_\PP(s_n^{2K})$ and $O_\PP(s_n^{K+2})$. For every $K\ge1$, their sum is $o_\PP(s_n^{K+1})$.

Subtracting \eqref{eq:propagator-error-estimator} and proceeding through the fixed number of subintervals gives
\[
\max_{j\le M}\|e_{j}^K-\widehat e_{j}^K\| =o_\PP(s_n^{K+1}).
\]
The partial fine maps have uniformly bounded second derivatives on the relevant compact sets. Their expansion around $U_{j}^{(K-1)}$ therefore gives, uniformly in $u$,
\[
\begin{aligned}
\Theta^{(K,M)}(T_j+u)-\Theta(T_j+u) &= D\mathcal F_{j,u}(U_{j}^{(K-1)})e_{j}^K \\
&\quad+\frac u\Delta \left\{ e_{j+1}^K - D\mathcal F_{j}(U_{j}^{(K-1)})e_{j}^K \right\} +o_\PP(s_n^{K+1}).
\end{aligned}
\]
Subtracting \eqref{eq:propagator-path-estimator} proves the assertion.
\end{proof}

If, in addition, $f_n$ and $\sigma_n$ have bounded third derivatives on $\Omega_n$, uniformly in $n$, the argument giving \eqref{eq:second-derivative-consistency} yields $\max_{j<M}\|D^2\delta_{j}\|_{C^0(D)}=O_\PP(s_n+\lambda+\tau)=O_\PP(s_n)$. The three remainder terms in the proof are then $O_\PP(s_n^{K+2})$, and
\[
\|\widetilde\Theta^{(K,M)}-\Theta\|_\infty =O_\PP(s_n^{K+2}).
\]

The estimator requires only directional derivatives. On each subinterval, $D\mathcal F_{j}(U_{j}^{(K-1)})v_{j}^K$ is obtained by propagating one tangent vector along the fine trajectory started at $U_{j}^{(K-1)}$, with the same learner draws; these sweeps run concurrently. The recursion \eqref{eq:propagator-error-estimator} is then a sequential sweep of coarse tangent updates, and \eqref{eq:propagator-path-estimator} needs a second concurrent fine tangent sweep along the same trajectories. In forward-mode differentiation \citep{GriewankWalther2008} each tangent sweep costs a fixed multiple of the corresponding propagation sweep, so the correction adds two concurrent fine sweeps and one sequential coarse sweep to the cost of level $K$.

\subsection{Examples}\label{subsec:explicit}

A scalar multiplicative learner gives an explicit law at every correction level. Let $a>0$ and $\gamma\in\R$, and consider
\begin{equation}\label{eq:critical-scalar}
\Theta_{\ell+1} =\Theta_{\ell} \{1-a\lambda+\gamma\lambda\bar\xi_{\ell+1}\}, \qquad \Theta_{0}=\beta_0,
\end{equation}
where the individual draws $\xi_{\ell,r}$ are independent Rademacher variables. The coarse propagator multiplier is
\[
g=(1-a\tau)^{J}.
\]
For sufficiently large $n$, $g>0$. Define
\[
F_{j} =\prod_{\ell=jL+1}^{(j+1)L} \{1-a\lambda+\gamma\lambda\bar\xi_{\ell}\}, \qquad z_{j}=\frac{F_{j}}{g}-1.
\]
The bounded draws and the fixed number of subintervals give uniform pathwise bounds for all states through level $M$. A smooth cutoff outside a larger ball therefore gives bounded coefficients without changing these trajectories.

For $0\le r\le M$, let $e_r$ be the elementary symmetric polynomial of degree $r$ in $M$ variables, with $e_0=1$.

\begin{theorem}\label{thm:scalar}
For the scalar model \eqref{eq:critical-scalar},
\begin{equation}\label{eq:scalar-exact}
\begin{aligned}
U_{M}^{(K)} &=\beta_0g^M \sum_{r=0}^{\min(K,M)} e_r(z_{0},\ldots,z_{M-1}),\\
\Theta(T) &=\beta_0g^M \prod_{j=0}^{M-1}(1+z_{j}).
\end{aligned}
\end{equation}
If $\tau/s_n\to c<\infty$, define
\[
Y_j=\gamma\{W(T_{j+1})-W(T_j)\} +\frac{ca^2\Delta}{2}, \qquad 0\le j<M.
\]
Then, for every $0\le K<M$,
\begin{equation}\label{eq:scalar-limit}
s_n^{-(K+1)}\{U_{M}^{(K)}-\Theta(T)\} \Rightarrow -\beta_0e^{-aT}e_{K+1}(Y_0,\ldots,Y_{M-1}).
\end{equation}
The convergences hold jointly over the correction levels and with $W_n$. For $K\ge M$, the error is identically zero.
\end{theorem}

\begin{proof}
The scalar Parareal recursion reads
\[
U_{j+1}^{(K+1)} =gU_{j}^{(K+1)}+gz_{j}U_{j}^{(K)}.
\]
Dividing by $g^{j+1}$ gives the recursion for partial sums of elementary symmetric polynomials. Induction proves \eqref{eq:scalar-exact}.

Since $\lambda/s_n\to0$, expansion of the products gives
\[
\frac{z_{j}}{s_n} =\gamma\{W_n(T_{j+1})-W_n(T_j)\} +\frac{\tau}{s_n}\frac{a^2\Delta}{2} +o_\PP(1).
\]
More precisely, expansion of the logarithms shows that the deterministic fine-step bias is $O(\lambda)$, the quadratic learner terms are $O_\PP(s_n^2)$, and the remaining coarse bias is $O(\tau^2)$. Each is negligible after division by $s_n$. The expansion holds jointly over the fixed collection of propagators, so
\[
(z_{j}/s_n)_{j=0}^{M-1} \Rightarrow(Y_j)_{j=0}^{M-1},
\]
jointly with $W_n$.

By \eqref{eq:scalar-exact},
\[
U_{M}^{(K)}-\Theta(T) =-\beta_0g^M \sum_{r=K+1}^{M} e_r(z_{0},\ldots,z_{M-1}).
\]
After normalization, the first term gives the limit and every higher-degree term vanishes. Finally, $g^M\to e^{-aT}$.
\end{proof}

When $c=0$, products indexed by distinct subsets of the independent centered Gaussian increments are orthogonal. The variance of the terminal limit is therefore
\[
\beta_0^2e^{-2aT} \binom{M}{K+1}(\gamma^2\Delta)^{K+1}.
\]
If $\beta_0\gamma\ne0$, this is positive for every $K<M$. Thus the exponent $K+1$ in the error bound cannot be improved in general.

For $c>0$, independence of the $Y_j$ gives the limiting mean
\[
-\beta_0e^{-aT} \binom{M}{K+1} \left(\frac{ca^2\Delta}{2}\right)^{K+1}.
\]
For one correction and $c=0$, the terminal limit is
\[
-\frac{\beta_0e^{-aT}\gamma^2}{2} \left\{ W(T)^2 -\sum_{j=0}^{M-1} \bigl(W(T_{j+1})-W(T_j)\bigr)^2 \right\}.
\]
For fixed $M$, the sum of squared propagator increments remains random.

We next specialize the theorem to the coordinate learner. Its additional sampling condition is collected in the following assumption.

\begin{assumption}\label{ass:coordinate-limit}
The coordinate probabilities $\pi_1,\ldots,\pi_p$ are fixed and strictly positive, and
\[
H_n\xrightarrow{\PP}H,\qquad b_n\xrightarrow{\PP}b
\]
for deterministic $H\in\R^{p\times p}$ and $b\in\R^p$.
\end{assumption}

For $p\ge2$, define
\[
v_\pi=(\sqrt{\pi_1},\ldots,\sqrt{\pi_p})^\top
\]
and choose a $p\times(p-1)$ matrix $Q$ satisfying
\[
Q^\top Q=I_{p-1},\qquad Q^\top v_\pi=0.
\]
For $1\le i\le p$, put
\[
D_i=\pi_i^{-1}e_ie_i^\top,\qquad d_{n,i}(x)=(D_i-I_p)(b_n-H_nx),\qquad \xi_i=\frac{Q^\top e_i}{\sqrt{\pi_i}}.
\]
For a coordinate draw $\iota$ with probabilities $\pi_i$,
\[
\E\xi_\iota=0,\qquad \E(\xi_\iota\xi_\iota^\top)=I_{p-1},
\]
and the innovations are bounded. Moreover,
\[
\widehat F(x,\iota)=f_n(x)+\sigma_n(x)\xi_\iota, \qquad f_n(x)=b_n-H_nx,
\]
where
\[
\sigma_n(x) = \bigl[ \sqrt{\pi_1}d_{n,1}(x),\ldots, \sqrt{\pi_p}d_{n,p}(x) \bigr]Q.
\]

\begin{corollary}\label{cor:coordinate-limit}
Under Assumption~\ref{ass:coordinate-limit}, with $p\ge2$ and $\tau/s_n\to c<\infty$, the joint limit \eqref{eq:multivariate-joint} and its conditional version hold for the coordinate learner. The limiting fields are
\[
f(x)=b-Hx,\qquad \sigma(x) = \bigl[ \sqrt{\pi_1}(D_1-I_p)f(x),\ldots, \sqrt{\pi_p}(D_p-I_p)f(x) \bigr]Q.
\]
No positive-definiteness assumption on $H$ is required.
\end{corollary}

\begin{proof}
Fix $R>\|H\|+\|b\|$ and define
\[
\Omega_{n,R}=\{\|H_n\|+\|b_n\|\le R\}.
\]
Its probability tends to one. The affine matrix-product bounds from Section~\ref{subsec:coordinate} give
\[
\E_\xi\|\mathcal F_{j}(x)\|^2 \le C_R(1+\|x\|^2), \qquad x\in\R^p,
\]
on $\Omega_{n,R}$. Thus the conditional integrability condition is satisfied.

The same product bounds control every intermediate fine and coarse state started in a bounded set. Induction through the finite collection of propagators and correction levels $K\le M$ places all states used by these levels in a deterministic ball. Its radius depends on $R$, $T$, $M$, $\beta_0$ and the coordinate probabilities, but not on $n$.

Choose a smooth compactly supported cutoff equal to one on a neighborhood of that ball. Multiplying $f_n$ and $\sigma_n$ by the cutoff gives $C_b^2$ fields for every training sample, with uniform bounds on $\Omega_{n,R}$. Assumption~\ref{ass:coordinate-limit} gives their convergence in $C^2$ on compact sets. The transformed coordinate draws have a fixed centered law, identity covariance and bounded support, and remain independent of the training data. Every condition of Assumption~\ref{ass:multi} is therefore satisfied by the cutoff fields.

The cutoff and original algorithms agree on $\Omega_{n,R}$. Their mean flows and the derivatives entering the limiting recursion also agree along the relevant trajectories. Theorem~\ref{thm:parareal-limit} therefore proves the result.
\end{proof}

For this affine drift,
\[
\Phi(t,u)=e^{-H(t-u)},\qquad \mathsf E_j=-\frac{\Delta}{2}e^{-H\Delta}H^2.
\]
Write $\sigma_r$ for the $r$-th column of $\sigma$ and set $S_r=D\sigma_r$, which is constant. Then
\[
\mathcal R_j[v] =\sum_{r=1}^{p-1} \int_{T_j}^{T_{j+1}} e^{-H(T_{j+1}-u)} S_r e^{-H(u-T_j)}v\,dW_r(u).
\]
These formulas determine the limiting recursion from $H$, $b$ and the coordinate probabilities. When $p=1$, the learner is deterministic and the path limits follow directly by setting $\sigma=0$.

\section{Statistical accuracy and computational cost}
\label{sec:statistics}

The error bounds give sufficient conditions for preserving the statistical limits of the sequential algorithm. We consider path estimation and validation-based selection, then calculate the propagation cost.

\subsection{Path limits}

Let $\theta_n$ and $\theta$ denote the empirical and population mean flows introduced in Section~\ref{sec:second}. Recall that
\[
\eta_n=s_n+\lambda+\tau.
\]

\begin{assumption}\label{ass:statistical-limit}
The learner satisfies the hypotheses of Proposition~\ref{prop:primitive} or Theorem~\ref{thm:parareal}. In the smooth setting, the fields remain $C_b^2$ almost surely, and the uniform coefficient bounds may instead hold on events $\Omega_n\in\mathcal D_n$ with $\PP(\Omega_n)\to1$, with constants independent of $n$.

There is a deterministic sequence $a_n\to\infty$ such that
\begin{equation}\label{eq:existing-limit}
\left( a_n(\theta_n-\theta),\, s_n^{-1}(\Theta-\theta_n) \right) \Rightarrow(U,V)
\end{equation}
in $C([0,T];\R^p)^2$, where $U$ and $V$ are independent continuous processes.
\end{assumption}

Theorem~2.6 of \citet{BuiEsstafaJoint2026} establishes the sequential joint limit under its empirical-field and learner assumptions. We take this statistical limit as input. The sampling process $U$ need not be Gaussian.

Write
\[
\mathcal U_n=a_n(\theta_n-\theta),\qquad r_n=a_ns_n,
\]
and, for each integer $K\ge0$, define
\[
\varepsilon_{n,K}=\eta_n^{K+1}.
\]
The estimates of Section~\ref{sec:bounds} give
\begin{equation}\label{eq:statistical-path-error}
\|\Theta^{(K,M)}-\Theta\|_\infty =O_\PP(\varepsilon_{n,K})
\end{equation}
for every fixed $K$.

\begin{corollary}\label{cor:CLT}
Under Assumption~\ref{ass:statistical-limit}, let $1\le K<M$ be fixed. If
\[
\varepsilon_{n,K}=o(s_n),
\]
then
\begin{equation}\label{eq:parareal-joint-statistical-limit}
\left( \mathcal U_n,\, s_n^{-1}(\Theta^{(K,M)}-\theta_n) \right) \Rightarrow(U,V).
\end{equation}
If $r_n\to\rho\in[0,\infty)$, then
\[
a_n(\Theta^{(K,M)}-\theta) \Rightarrow U+\rho V.
\]
If $r_n\to\infty$, then
\[
s_n^{-1}(\Theta^{(K,M)}-\theta) \Rightarrow V.
\]
All convergences hold in continuous-path spaces. For $K\ge M$, these conclusions hold without a numerical rate condition.
\end{corollary}

\begin{proof}
For $K<M$, equation~\eqref{eq:statistical-path-error} gives
\[
s_n^{-1} \|\Theta^{(K,M)}-\Theta\|_\infty \xrightarrow{\PP}0.
\]
Slutsky's theorem proves \eqref{eq:parareal-joint-statistical-limit}. When $r_n$ has a finite limit, use
\[
a_n(\Theta^{(K,M)}-\theta) =\mathcal U_n +r_n s_n^{-1} (\Theta^{(K,M)}-\theta_n).
\]
When $r_n\to\infty$, use
\[
s_n^{-1}(\Theta^{(K,M)}-\theta) =r_n^{-1}\mathcal U_n +s_n^{-1}(\Theta^{(K,M)}-\theta_n).
\]
For $K\ge M$, the paths coincide by Proposition~\ref{prop:finite-exactness}.
\end{proof}

\subsection{Validation-based stopping}
\label{subsec:stopping}

Fix a closed interval $I=[t_-,t_+]\subset(0,T)$. Let $R_v$ be the population validation loss and $R_{v,n}$ its empirical counterpart. The same empirical loss is used for both algorithms. Define
\[
C(t)=R_v(\theta(t)),\qquad C_n(t)=R_{v,n}(\Theta(t)),\qquad C_n^{(K,M)}(t) =R_{v,n}(\Theta^{(K,M)}(t)).
\]
For a real-valued function $g$ on $I$, write
\[
\|g\|_{\infty,I}=\sup_{t\in I}|g(t)|.
\]
Define
\[
\zeta_n=a_n^{3/2}s_n,\qquad h_n=a_n^{-1}+s_n^{2/3}.
\]

\begin{assumption}\label{ass:stopping}
The function $C$ is continuous on $I$ and has a unique global minimizer $t_\star\in\operatorname{int}(I)$. It is twice continuously differentiable near $t_\star$, with
\[
c_\star=C''(t_\star)>0.
\]

Each $R_{v,n}$ is continuous. There are a fixed closed ball $D_v$, whose interior contains $\theta(I)$, and nonnegative random variables $L_{v,n}=O_\PP(1)$ such that
\[
|R_{v,n}(x)-R_{v,n}(y)| \le L_{v,n}\|x-y\|, \qquad x,y\in D_v.
\]

The sequence $\zeta_n$ converges to $\zeta\in[0,\infty]$. Define
\[
\bar h_n=
\begin{cases}
a_n^{-1},&\zeta<\infty,\\
s_n^{2/3},&\zeta=\infty.
\end{cases}
\]
There are a deterministic constant $\sigma_\star>0$ and a standard two-sided Brownian motion $W_2$ with independent positive and negative halves. In the finite-$\zeta$ case, there is also a real-valued random variable $G_\star$, and $W_2$ is independent of $(U,G_\star)$. When $\zeta=\infty$, $W_2$ is independent of $U$. Define, for $u\in\R$,
\[
\mathcal Q_\zeta(u)=
\begin{cases}
\dfrac12c_\star u^2+uG_\star +\zeta\sigma_\star W_2(u), &\zeta<\infty,\\[1mm] \dfrac12c_\star u^2+\sigma_\star W_2(u), &\zeta=\infty.
\end{cases}
\]

For $u$ such that $t_\star+\bar h_nu\in I$, set
\[
\mathcal Q_n(u) =\bar h_n^{-2} \{C_n(t_\star+\bar h_nu)-C_n(t_\star)\}.
\]
For every $A<\infty$,
\begin{equation}\label{eq:sequential-local-criterion}
\left(\mathcal U_n,\mathcal Q_n|_{[-A,A]}\right) \Rightarrow \left(U,\mathcal Q_\zeta|_{[-A,A]}\right)
\end{equation}
in $C([0,T];\R^p)\times C([-A,A])$.

Finally, for every $b>0$,
\begin{equation}\label{eq:outer}
\lim_{R\to\infty}\limsup_{n\to\infty} \PP\left( \inf_{\substack{t\in I\\|t-t_\star|>R\bar h_n}} \{C_n(t)-C_n(t_\star)\} \le b\bar h_n^2 \right)=0.
\end{equation}
The infimum over an empty set is $+\infty$.
\end{assumption}

The local limits in \eqref{eq:sequential-local-criterion} are those of the sequential validation theory in \citet{BuiEsstafaJoint2026}. The random slope $G_\star$ may depend on $U$. Condition~\eqref{eq:outer} is imposed separately. It controls distant near-minimizers and is used to localize both validation selections.

Let $\widehat t_n$ and $\widehat t_n^{(K,M)}$ be the smallest minimizers of $C_n$ and $C_n^{(K,M)}$ over $I$. They exist because the criteria are continuous. Define the corresponding fitted values by
\[
\widehat\theta_n=\Theta(\widehat t_n),\qquad \widehat\theta_n^{(K,M)} =\Theta^{(K,M)}(\widehat t_n^{(K,M)}).
\]
For $K\ge M$, the paths and criteria coincide. The common tie-breaking rule then gives
\[
\widehat t_n^{(K,M)}=\widehat t_n,\qquad \widehat\theta_n^{(K,M)}=\widehat\theta_n.
\]

The two local scales satisfy
\begin{equation}\label{eq:stopping-scale-comparison}
\frac{h_n}{\bar h_n}\longrightarrow
\begin{cases}
1+\zeta^{2/3},&\zeta<\infty,\\
1,&\zeta=\infty.
\end{cases}
\end{equation}

\begin{proposition}\label{prop:risk}
Under Assumptions~\ref{ass:statistical-limit} and~\ref{ass:stopping}, for every fixed $K$,
\begin{equation}\label{eq:criterion-perturb}
\|C_n^{(K,M)}-C_n\|_{\infty,I} =O_\PP(\varepsilon_{n,K}).
\end{equation}
If $\varepsilon_{n,K}=o(h_n^2)$, then, for every $A<\infty$,
\[
\sup_{|u|\le A} |\mathcal Q_n^{(K,M)}(u)-\mathcal Q_n(u)| \xrightarrow{\PP}0,
\]
where
\[
\mathcal Q_n^{(K,M)}(u) =\bar h_n^{-2} \{C_n^{(K,M)}(t_\star+\bar h_nu) -C_n^{(K,M)}(t_\star)\}.
\]
\end{proposition}

\begin{proof}
Assumption~\ref{ass:statistical-limit} implies $\|\Theta-\theta\|_\infty=o_\PP(1)$. Equation~\eqref{eq:statistical-path-error} gives the same consistency for the Parareal path. Both paths therefore lie in $D_v$ over $I$ with probability tending to one. On that event,
\[
\|C_n^{(K,M)}-C_n\|_{\infty,I} \le L_{v,n} \|\Theta^{(K,M)}-\Theta\|_\infty.
\]
This proves \eqref{eq:criterion-perturb}. The centered local difference is bounded by
\[
2\bar h_n^{-2} \|C_n^{(K,M)}-C_n\|_{\infty,I},
\]
which tends to zero by \eqref{eq:stopping-scale-comparison}.
\end{proof}

Let
\[
T_\zeta=\argmin_{u\in\R}\mathcal Q_\zeta(u).
\]
The quadratic term dominates the Brownian and linear terms at infinity, so the minimum is attained. For $\zeta=0$, its unique location is $T_0=-G_\star/c_\star$. In the Brownian cases, almost sure uniqueness follows from Lemma~2.6 of \citet{KimPollard1990}, conditioning on $G_\star$ when necessary.

\begin{theorem}\label{thm:stopping}
Under Assumptions~\ref{ass:statistical-limit} and~\ref{ass:stopping}, let $0\le K<M$ be fixed and suppose
\[
\varepsilon_{n,K}=o(h_n^2).
\]
Then
\begin{equation}\label{eq:selected-time-equivalence}
\frac{\widehat t_n^{(K,M)}-\widehat t_n}{\bar h_n} \xrightarrow{\PP}0.
\end{equation}
If $\zeta<\infty$,
\begin{equation}\label{eq:selected-fit-finite}
\left( a_n(\widehat t_n^{(K,M)}-t_\star),\, a_n\{\widehat\theta_n^{(K,M)}-\theta(t_\star)\} \right) \Rightarrow \left( T_\zeta,\, U(t_\star)+\dot\theta(t_\star)T_\zeta \right).
\end{equation}
If $\zeta=\infty$,
\begin{equation}\label{eq:selected-fit-infinite}
\left( s_n^{-2/3}(\widehat t_n^{(K,M)}-t_\star),\, s_n^{-2/3} \{\widehat\theta_n^{(K,M)}-\theta(t_\star)\} \right) \Rightarrow \left( T_\infty,\, \dot\theta(t_\star)T_\infty \right).
\end{equation}
Moreover,
\begin{equation}\label{eq:selected-fit-equivalence}
\frac{ \|\widehat\theta_n^{(K,M)}-\widehat\theta_n\| }{\bar h_n} \xrightarrow{\PP}0.
\end{equation}
For $K\ge M$, the selections and fitted values agree exactly, so no numerical rate condition is needed.
\end{theorem}

\begin{proof}
Write
\[
\varpi_n=\|C_n^{(K,M)}-C_n\|_{\infty,I}.
\]
Proposition~\ref{prop:risk} and \eqref{eq:stopping-scale-comparison} give $\varpi_n=o_\PP(\bar h_n^2)$. By the minimizing property,
\[
C_n(\widehat t_n^{(K,M)})-C_n(t_\star) \le2\varpi_n.
\]
Consequently, for $R,b>0$,
\[
\begin{aligned}
&\PP\bigl( |\widehat t_n^{(K,M)}-t_\star|>R\bar h_n \bigr)\\
&\quad\le \PP\left( \inf_{\substack{t\in I\\|t-t_\star|>R\bar h_n}} \{C_n(t)-C_n(t_\star)\} \le b\bar h_n^2 \right) +\PP(2\varpi_n>b\bar h_n^2).
\end{aligned}
\]
The separation condition makes $(\widehat t_n^{(K,M)}-t_\star)/\bar h_n$ tight. The same argument, with $C_n(\widehat t_n)\le C_n(t_\star)$, makes $(\widehat t_n-t_\star)/\bar h_n$ tight.

The local criteria converge jointly as
\[
(\mathcal U_n,\mathcal Q_n,\mathcal Q_n^{(K,M)}) \Rightarrow(U,\mathcal Q_\zeta,\mathcal Q_\zeta)
\]
uniformly on compact local-time intervals. Every fixed compact interval eventually lies in the rescaled selection domain because $t_\star$ is interior to $I$. Tightness and uniqueness of the limiting minimizer identify the joint limit as
\[
\left( \mathcal U_n,\, \frac{\widehat t_n-t_\star}{\bar h_n},\, \frac{\widehat t_n^{(K,M)}-t_\star}{\bar h_n} \right) \Rightarrow(U,T_\zeta,T_\zeta).
\]
This is the argmin argument of Theorem~2.7 in \citet{KimPollard1990}, applied to the two criteria with their common limit. It proves \eqref{eq:selected-time-equivalence}.

The joint path limit gives
\[
\|\Theta-\theta_n\|_\infty=O_\PP(s_n),
\]
and $\mathcal U_n$ is stochastically equicontinuous. The numerical error is $o_\PP(\bar h_n^2)$ and hence $o_\PP(\bar h_n)$.

If $\zeta<\infty$, then
\[
a_ns_n=\zeta_n a_n^{-1/2}\longrightarrow0.
\]
Uniformly for $|u|\le A$, differentiability of $\theta$ and stochastic equicontinuity give
\[
\begin{aligned}
a_n\{ \Theta^{(K,M)}(t_\star+a_n^{-1}u) -\theta(t_\star) \} =\mathcal U_n(t_\star) +\dot\theta(t_\star)u+o_\PP(1).
\end{aligned}
\]
The sequential path has the same expansion. Evaluating at the two tight normalized selections proves \eqref{eq:selected-fit-finite} and \eqref{eq:selected-fit-equivalence}.

If $\zeta=\infty$, then
\[
\frac{a_n^{-1}}{s_n^{2/3}} =\zeta_n^{-2/3}\longrightarrow0, \qquad \frac{s_n}{s_n^{2/3}} =s_n^{1/3}\longrightarrow0.
\]
Thus, uniformly for bounded $u$,
\[
s_n^{-2/3} \{ \Theta^{(K,M)}(t_\star+s_n^{2/3}u) -\theta(t_\star) \} =\dot\theta(t_\star)u+o_\PP(1).
\]
Again the sequential path has the same expansion. The joint convergence of the normalized selections proves \eqref{eq:selected-fit-infinite} and \eqref{eq:selected-fit-equivalence}.

For $K\ge M$, finite-step exactness gives the remaining assertion.
\end{proof}

\begin{corollary}\label{cor:one-statistical}
Under Assumptions~\ref{ass:statistical-limit} and~\ref{ass:stopping}, suppose $\lambda/\tau\to0$ and $\tau=O(s_n)$. Then every fixed $K\ge1$ preserves both the joint path limit of Corollary~\ref{cor:CLT} and the validation-selection limits of Theorem~\ref{thm:stopping}. In particular, one Parareal correction suffices.
\end{corollary}

\begin{proof}
For $K\ge M$, use finite-step exactness. For $1\le K<M$, the two mesh conditions give $\lambda=(\lambda/\tau)(\tau/s_n)s_n=o(s_n)$. Hence $\eta_n=O(s_n)$ and
\[
\frac{\varepsilon_{n,K}}{s_n} =O(s_n^K)\longrightarrow0.
\]
Since $h_n\ge s_n^{2/3}$,
\[
\frac{\varepsilon_{n,K}}{h_n^2} =O(s_n^{K-1/3})\longrightarrow0.
\]
The two preceding results apply.
\end{proof}

\subsection{Computational cost}\label{sec:computational}

We compare propagation costs after preprocessing shared by both algorithms. Validation-loss evaluations are excluded from both costs.

Fix $0\le K<M$. The model has $M$ processors, one for each subinterval. A coarse sweep is sequential, while the fine propagators in a sweep run concurrently. Coarse values from the preceding correction level are cached.

Let $c_{F,n}>0$ and $c_{G,n}>0$ be the costs of one fine update and one coarse Euler update. A fine propagator costs $L c_{F,n}$ and a coarse propagator costs $J c_{G,n}$. Let $c_{{\rm pre},n}\ge0$ collect additional setup and output costs, and let $c_{{\rm comm},n}\ge0$ be the cost of one correction's endpoint exchange and synchronization.

We count full sweeps without overlap between coarse and fine work. There are $K+1$ coarse sweeps, $K$ fine sweeps for the corrections, and one additional fine sweep for the reconstruction. The critical-path cost is
\begin{equation}\label{eq:cost}
\mathcal T_{{\rm par},n}(M,K) =c_{{\rm pre},n} +(K+1)M J c_{G,n} +(K+1)L c_{F,n} +Kc_{{\rm comm},n}.
\end{equation}
The sequential propagation cost is
\[
\mathcal T_{{\rm seq},n} =N c_{F,n}=M L c_{F,n}.
\]
The model does not include load imbalance or memory effects.

Define the coarse-to-fine propagator cost ratio
\[
\chi_n =\frac{J c_{G,n}}{L c_{F,n}} =\frac{\lambda}{\tau} \frac{c_{G,n}}{c_{F,n}}.
\]
Let $\mathcal S_n(M,K)$ be the ratio of the sequential cost to \eqref{eq:cost}. Let $\mathcal S_n^0(M,K)$ denote the same ratio when $c_{{\rm pre},n}=c_{{\rm comm},n}=0$.

\begin{proposition}\label{prop:cost}
In the propagation model \eqref{eq:cost},
\begin{equation}\label{eq:fixed-M-speedup}
\mathcal S_n^0(M,K) =\frac{M}{(K+1)(1+M\chi_n)},
\end{equation}
and
\[
\mathcal S_n(M,K) \le\mathcal S_n^0(M,K) \le\frac{M}{K+1}.
\]
If
\[
\chi_n\longrightarrow0,\qquad c_{{\rm pre},n}+Kc_{{\rm comm},n} =o(L c_{F,n}),
\]
then
\[
\mathcal S_n(M,K)\longrightarrow\frac{M}{K+1}.
\]
\end{proposition}

\begin{proof}
Dividing the costs by $L c_{F,n}$ gives
\[
\mathcal S_n(M,K) = \frac{M}{ (K+1)(1+M\chi_n) +\{c_{{\rm pre},n}+Kc_{{\rm comm},n}\}/ (L c_{F,n}) }.
\]
The formula, bounds and limit follow.
\end{proof}

The model corresponds to a reconstruction sweep in which each processor stores its own fine path, as in the timing experiment of Section~\ref{subsec:num-performance}. Validation can be added in the same way: if $c_{V,n}$ is the cost of evaluating the validation criterion along one fine step, the sequential cost becomes $N(c_{F,n}+c_{V,n})$ and each processor adds $L c_{V,n}$ to its reconstruction sweep.

The factor $K+1$ includes the reconstruction sweep. With $M=6$, the propagation model therefore bounds the speed-up by $3$, $2$ and $3/2$ for $K=1,2,3$. If only subinterval endpoints are required, the reconstruction sweep is omitted and the fine cost becomes $K L c_{F,n}$.

For gradients evaluated observation by observation, the cost ratio can be expressed in training epochs.

\begin{corollary}\label{cor:epochs}
Suppose a fine update averages $B$ per-observation gradients and a coarse update evaluates the full gradient over $n$ observations. Taking one per-observation gradient as the cost unit gives
\[
c_{F,n}=B,\qquad c_{G,n}=n.
\]
Define the number of sequential training epochs by
\[
E_{n,\mathrm{epoch}}=\frac{NB}{n}.
\]
Then
\begin{equation}\label{eq:epoch-speedup}
\chi_n=\frac{n\lambda}{B\tau} =\frac{M J}{E_{n,\mathrm{epoch}}}, \qquad \mathcal S_n^0(M,K) =\frac{E_{n,\mathrm{epoch}}}{ (K+1)(M J+E_{n,\mathrm{epoch}}/M) }.
\end{equation}
If $M\le K+1$, the ideal speedup is less than one. For $M>K+1$, it exceeds one exactly when
\[
E_{n,\mathrm{epoch}}> \frac{(K+1)M^2J}{M-K-1}.
\]
\end{corollary}

\begin{proof}
The identities $N=M L$ and $J/L=\lambda/\tau$ give the expression for $\chi_n$. Substitution into \eqref{eq:fixed-M-speedup} gives \eqref{eq:epoch-speedup}. The speedup exceeds one precisely when
\[
E_{n,\mathrm{epoch}}\left(1-\frac{K+1}{M}\right) >(K+1)M J.
\]
This is impossible for $M\le K+1$. For $M>K+1$, division by the positive coefficient gives the stated threshold.
\end{proof}

One coarse sweep uses $M J=T/\tau$ full-gradient passes. Its cost must therefore be compared with $E_{n,\mathrm{epoch}}$. Accuracy and cost are linked through $\tau$. The conditions $\eta_n^{K+1}=o(s_n)$ and $\eta_n^{K+1}=o(h_n^2)$ are sufficient to preserve the joint learner limit and the validation-selection limit, respectively, while the cost ratio contains $\lambda/\tau$.

For the coordinate learner, the data matrices $H_n$ and $b_n$ can be formed once and shared by both algorithms. With dense $H_n$, a direct implementation uses $O(Bp)$ operations per fine update and $O(p^2)$ operations per coarse update. For fixed $p$ and $B$, the per-update costs remain bounded above and away from zero. Consequently, $J/L\to0$ gives $\chi_n\to0$. Under the overhead condition of Proposition~\ref{prop:cost}, the propagation speedup approaches $M/(K+1)$.

\section{Numerical experiments}\label{sec:numerical}

We examine Parareal errors for two sample-dependent $L^2$ residual-fitting learners, compare their normalized errors with the second-order limit law, and measure the elapsed time of a parallel implementation. In the convergence experiments, both numerical discretizations are coupled to the training-sample size $n$. The timing experiment solves the largest of these problems once.

\subsection{Model, discretization and error measures}\label{subsec:num-design}

All experiments use $p=4$ covariates, a time horizon $T=4$, $M=6$ subintervals, batch size $B=1$, and $\beta_0=0$. Independent training and validation samples both have size $n$. The training observations satisfy
\begin{equation}\label{eq:num-regression-model}
X_i\sim N(0,\Sigma),\qquad Y_i=1.6X_i^\top\beta_v+\varepsilon_i, \qquad \Sigma_{rs}=0.5^{|r-s|},
\end{equation}
where $\beta_v=(1,-0.8,0.6,0.4)^\top$ and the $\varepsilon_i$ are independent standard Gaussian variables. Validation observations have the same covariate and noise distributions but regression coefficient $\beta_v$. Write
\[
H_n=\frac1n\sum_{i=1}^n X_iX_i^\top,\qquad b_n=\frac1n\sum_{i=1}^n X_iY_i,
\]
with $H_{v,n}$ and $b_{v,n}$ defined from the independent validation sample. Up to an additive constant, the validation criterion is
\[
R_{v,n}(x)=\tfrac12x^\top H_{v,n}x-b_{v,n}^\top x, \qquad I=[0.4,3.6].
\]
The distinct training and validation coefficients yield an interior minimum of the population validation criterion for each of the two mean flows.

For the \emph{coordinate learner}, draw $i$ uniformly from $\{1,\ldots,p\}$ and set
\begin{equation}\label{eq:num-coordinate-learner}
\widehat F(x,i) =p\{b_{n,i}-(H_nx)_i\}e_i, \qquad f_n(x)=b_n-H_nx.
\end{equation}
For the \emph{nonlinear residual-fitting learner}, let $\alpha_{n,i}(x)$ be the unique minimizer over $a\in\mathbb R$ of
\[
\frac1{2n}\|Y-Xx-aX_{\cdot i}\|^2 +\frac\rho2 a^2+\frac\gamma4 a^4, \qquad \rho=0.25,\quad\gamma=0.5.
\]
Equivalently, $\alpha=\alpha_{n,i}(x)$ solves
\begin{equation}\label{eq:num-nonlinear-fit}
(H_{n,ii}+\rho)\alpha+\gamma\alpha^3 =b_{n,i}-(H_nx)_i.
\end{equation}
The update and its conditional mean are
\begin{equation}\label{eq:num-nonlinear-learner}
\widehat F(x,i)=p\alpha_{n,i}(x)e_i, \qquad f_n(x)=(\alpha_{n,1}(x),\ldots,\alpha_{n,p}(x))^\top.
\end{equation}
The penalty makes the fit unique and the update map smooth and nonlinear in the current coefficient vector $x$; the selected predictor $aX_{\cdot i}$ remains linear in its covariate. Both learners are built from the current training sample. Their untruncated fields have linear growth rather than the global boundedness of Assumption~\ref{ass:multi}; the localization argument is recorded in Appendix~\ref{app:num-localization}.

We choose the sample sizes and the two grids jointly. With $\Delta=T/M=2/3$,
\begin{equation}\label{eq:num-joint-grid}
\begin{gathered}
n\in\{12^2,15^2,18^2,24^2,30^2,36^2,42^2\},\\
N=4n^{3/2},\quad L=N/M,\quad \lambda=T/N=n^{-3/2},\quad s_n=\sqrt{\lambda}=n^{-3/4},\\
J=\operatorname{round}(\Delta/s_n),\qquad \tau=\Delta/J.
\end{gathered}
\end{equation}
Thus $n$ increases from $144$ to $1764$, $\tau/s_n\to1$ and $\lambda/\tau\to0$. With $a_n=\sqrt n$,
\[
a_n^{3/2}s_n=1,\qquad h_n=n^{-1/2}+s_n^{2/3}=2n^{-1/2},\qquad \delta_n=n^{-1/2}+s_n.
\]

The fine learner draws are reused at every Parareal correction. Both learners use the same training--validation samples and coordinate-index sequences at matching sample sizes and replication indices. We use $R=1000$ independent replications for each $n$ and report $K=1,2,3$ corrections. Let $\Theta$ be the sequential fine path and $\Theta^{(K,M)}$ the continuously reconstructed Parareal path. The path and validation-criterion errors are
\begin{equation}\label{eq:num-errors}
E_n^K=\|\Theta^{(K,M)}-\Theta\|_\infty, \qquad D_n^K=\|C_n^{(K,M)}-C_n\|_{\infty,I}.
\end{equation}
We use the Euclidean norm in $\mathbb R^p$. Path suprema are computed on the fine grid for piecewise-affine reconstructions; criterion suprema also check stationary points inside each segment. The selected validation time minimizes the continuous interpolated criterion over $I$, taking the earliest minimizer in the case of ties. For any nonnegative error $Y$ we report the empirical root mean square $\{R^{-1}\sum_{r=1}^R Y_r^2\}^{1/2}$.

\subsection{Convergence and statistical accuracy}\label{subsec:num-statistics}

\begin{figure}[htbp]
\centering
\includegraphics[width=0.9\linewidth]{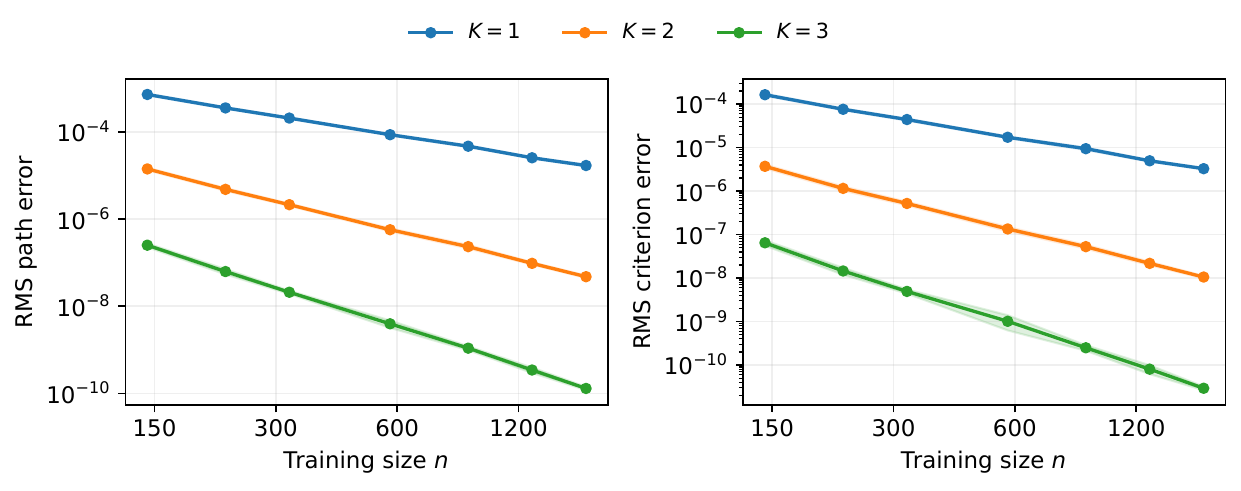}
\caption{Coordinate learner: RMS path error (left) and validation-criterion error (right) against training-sample size $n$, for $K=1,2,3$, with $1000$ independent replications per size. Axes are logarithmic; shaded bands are pointwise approximate 95\% Monte Carlo intervals.} \label{fig:num-coordinate}
\end{figure}

\begin{figure}[htbp]
\centering
\includegraphics[width=0.9\linewidth]{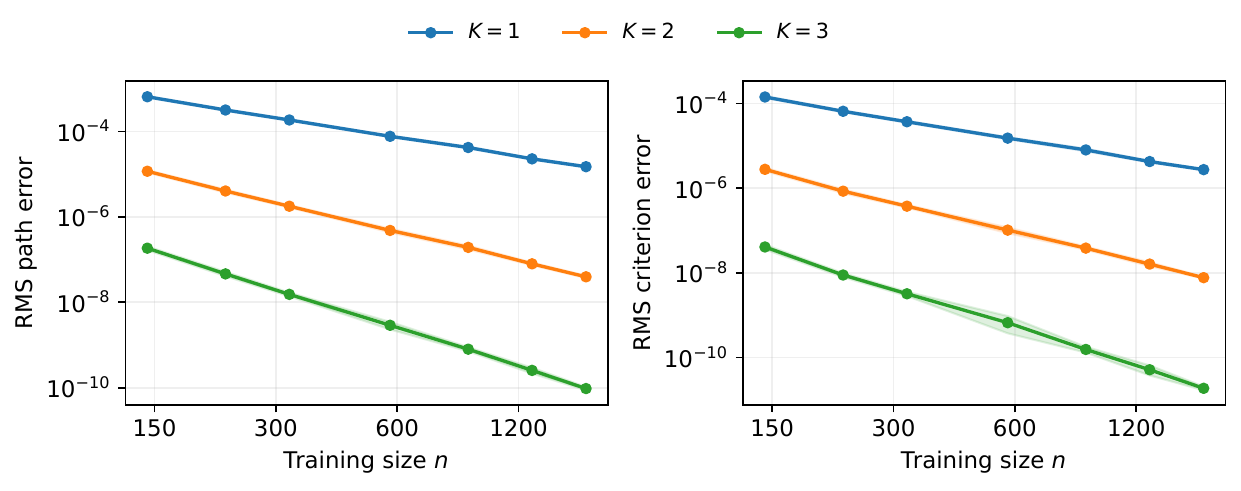}
\caption{Nonlinear residual-fitting learner, with the same sample sizes, replications and conventions as in Figure~\ref{fig:num-coordinate}.} \label{fig:num-nonlinear}
\end{figure}

In the critical regime $\tau/s_n\to1$, the theoretical path-error scale after $K$ corrections is $s_n^{K+1}=n^{-3(K+1)/4}$. Figures~\ref{fig:num-coordinate} and~\ref{fig:num-nonlinear} plot the RMS path errors $E_n^K$ and validation errors $D_n^K$ against $n$. Shaded bands are pointwise approximate 95\% Monte Carlo intervals based on delta-method standard errors.

Table~\ref{tab:num-joint-rates} gives the slopes obtained by regressing the seven log RMS values on $\log n$. The path slopes agree with the reference exponents $-3(K+1)/4$ to two decimals for both learners, and the normalized RMS $\mathrm{RMS}(E_n^K/s_n^{K+1})$ remains of order one at $n=1764$. The criterion slopes are slightly steeper than the path slopes. At $n=1764$ the RMS of $E_n^K/\delta_n$ is $6.2\times10^{-4}$ for the coordinate learner and $5.4\times10^{-4}$ for the nonlinear learner when $K=1$, and the RMS of $D_n^K/h_n^2$ is $1.4\times10^{-3}$ and $1.2\times10^{-3}$. One correction is therefore already far below the learner and stopping scales.

\begin{table}[htbp]
\centering
\small
\setlength{\tabcolsep}{4pt}
\caption{Path and criterion error slopes fitted from seven log--log RMS points, with $1000$ independent replications per sample size. The sixth column is the RMS of $E_n^K/s_n^{K+1}$ at $n=1764$; the last column counts changed selected times over all $7\times1000$ replications.} \label{tab:num-joint-rates}
\begin{tabular}{lrrrrrr}
\toprule
Learner & $K$ & \shortstack{Reference\\exponent} & \shortstack{Path\\slope} & \shortstack{Criterion\\slope} & \shortstack{RMS of\\$E_n^K/s_n^{K+1}$} & \shortstack{Changed\\times} \\
\midrule
\multirow{3}{*}{Coordinate} & 1 & $-1.50$ & $-1.50$ & $-1.56$ & $1.268$ & 40 \\
& 2 & $-2.25$ & $-2.25$ & $-2.31$ & $0.965$ & 2 \\
& 3 & $-3.00$ & $-3.00$ & $-3.03$ & $0.712$ & 0 \\
\addlinespace[3pt]
\multirow{3}{*}{Nonlinear} & 1 & $-1.50$ & $-1.50$ & $-1.56$ & $1.108$ & 25 \\
& 2 & $-2.25$ & $-2.25$ & $-2.31$ & $0.797$ & 2 \\
& 3 & $-3.00$ & $-2.99$ & $-3.02$ & $0.525$ & 1 \\
\bottomrule
\end{tabular}
\end{table}

For the coordinate and nonlinear population mean flows, the validation minimizers on $I$ are $1.9280$ and $2.5636$, with second derivatives of the population criterion $0.0454$ and $0.0273$. A selected time is counted as changed when $|\widehat t_n^{(K,M)}-\widehat t_n|>10^{-8}h_n$. For $K=1$, there are $40$ changed times for the coordinate learner and $25$ for the nonlinear learner, most of them at the smallest sample sizes. These occur when the empirical criterion has two nearly equal local minima; the largest changes are $5.9\,h_n$ and $3.8\,h_n$, respectively. For $K=2$ there are $2$ and $2$ changed times, and for $K=3$ there are $0$ and $1$.

\subsection{Second-order error distributions}\label{subsec:num-distributions}

We also examine the first coordinate of the signed terminal error after one correction,
\begin{equation}\label{eq:num-signed-limit}
Y_n=s_n^{-2}\bigl[\Theta^{(1,M)}(T)-\Theta(T)\bigr]_1 \ \Rightarrow\ [Z^{1,1}(T)]_1,
\end{equation}
where the limit is given by Theorem~\ref{thm:parareal-limit} and \eqref{eq:general-limit-recursion}. For each learner, we simulate $50\,000$ independent limit values using the \emph{population} coefficient fields and the joint Gaussian propagator increments and sensitivities from \eqref{eq:propagator-gaussian-increment}--%
\eqref{eq:propagator-random-sensitivity}. No parameter of this reference is fitted to the finite-$n$ errors. Details are in Appendix~\ref{app:num-limit-reference}.

\begin{figure}[htbp]
\centering
\includegraphics[width=0.9\linewidth]{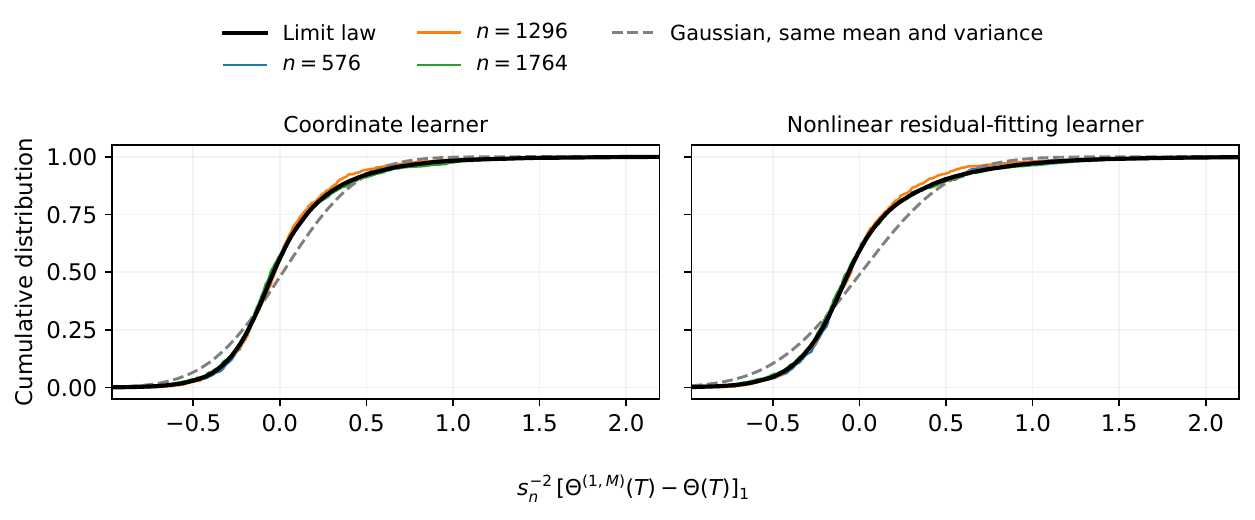}
\caption{Second-order terminal-error distributions for the coordinate (left) and nonlinear residual-fitting (right) learners. Empirical CDFs of $Y_n$ at $n=576,1296,1764$ use $1000$ independent replications per size. The limit CDF uses $50\,000$ independent draws of $[Z^{1,1}(T)]_1$. The dashed Gaussian CDF has the mean and variance of the simulated limit.} \label{fig:num-second-order-cdf}
\end{figure}

Figure~\ref{fig:num-second-order-cdf} shows agreement with the simulated second-order laws at all three sample sizes. At $n=1764$, the Kolmogorov distances between the empirical and limit CDFs are $0.033$ for the coordinate learner and $0.028$ for the nonlinear learner, below the 95\% Dvoretzky--Kiefer--Wolfowitz bound $0.043$ for $1000$ draws. The limit laws have skewness $1.60$ and $1.74$, and differ visibly from the Gaussian laws with the same mean and variance.

\FloatBarrier
\subsection{Computational performance}\label{subsec:num-performance}

The timing experiment uses the largest configuration of the convergence study, $n=1764$, $N=296\,352$, $M=6$ and $J=181$, with the data and learner draws of the first replication. The sequential reference runs the fine recursion once and minimizes the validation criterion along the way, without storing the path. Parareal runs on $P\in\{1,2,3,6\}$ worker processes, each responsible for $M/P$ consecutive subintervals. The coarse sweeps run on the main process. Each correction sends the current starting states to the workers, which return the fine endpoints. In the final sweep, each worker recomputes its fine paths from the corrected starting states, stores them locally, applies the endpoint matching \eqref{eq:bridge}, and minimizes the validation criterion over its subintervals; the main process keeps the smallest value. Each subinterval thus receives $(K+1)L$ fine updates, as in the model \eqref{eq:cost}.

Each configuration was timed $15$ times in randomized order after $3$ warm-up runs, and we report medians. Process creation and compilation are excluded. Numerical libraries are restricted to one thread, and each process is pinned to one physical core of an Intel Core i7-13800H laptop (six performance cores) running Ubuntu 24.04, with Python 3.12.3, NumPy 2.5.4 and Numba 0.68.1. Every timed run returned the sequential selected time to within $10^{-12}$.

\begin{table}[htbp]
\centering
\small
\setlength{\tabcolsep}{4pt}
\caption{Elapsed times in milliseconds, including validation selection, with $P=M=6$ processes, and errors of the same Parareal paths.} \label{tab:num-timing}
\begin{tabular}{lrrrrrr}
\toprule
Learner & $K$ & Sequential & Parareal & Speedup & $E_n^K/\delta_n$ & $D_n^K/h_n^2$ \\
\midrule
\multirow{3}{*}{Coordinate} & 1 & \multirow{3}{*}{$3.57$} & $2.51$ & $1.42$ & $3.3\times10^{-4}$ & $1.2\times10^{-3}$ \\
& 2 & & $3.25$ & $1.10$ & $1.1\times10^{-6}$ & $4.2\times10^{-6}$ \\
& 3 & & $3.71$ & $0.96$ & $2.9\times10^{-9}$ & $8.6\times10^{-9}$ \\
\addlinespace[3pt]
\multirow{3}{*}{Nonlinear} & 1 & \multirow{3}{*}{$17.40$} & $8.01$ & $2.17$ & $3.7\times10^{-4}$ & $1.1\times10^{-3}$ \\
& 2 & & $11.63$ & $1.50$ & $6.5\times10^{-7}$ & $1.9\times10^{-6}$ \\
& 3 & & $15.14$ & $1.15$ & $1.6\times10^{-9}$ & $3.7\times10^{-9}$ \\
\bottomrule
\end{tabular}
\end{table}

Table~\ref{tab:num-timing} reports the six-process results. For the nonlinear learner, one correction gives a speed-up of $2.17$, compared with the bound $M/(K+1)=3$ of Proposition~\ref{prop:cost}, while the path and criterion errors are far below $\delta_n$ and $h_n^2$. A coordinate update costs about a tenth of a nonlinear one, so the sequential run takes only $3.57$ ms and communication between processes takes a larger share of the parallel time; the speed-up for $K=1$ is $1.42$. Results for all processor counts, and the speed-ups predicted from measured unit costs, are in Appendix~\ref{app:num-timing}.

\appendix
\section{Additional numerical details}\label{app:numerics}

\subsection{Localization}\label{app:num-localization}

The training and validation covariates in Section~\ref{subsec:num-design} are Gaussian. The empirical coordinate and nonlinear residual-fit fields, therefore, are not uniformly bounded as functions of the state. We record the localization used to interpret the experiments through Assumption~\ref{ass:multi}. On the events
\[
\Omega_{n,R}=\{\|H_n\|+\|b_n\|\le R\},
\]
with $R$ sufficiently large, the probabilities tend to one by the law of large numbers. The coordinate update is affine with coefficients bounded on $\Omega_{n,R}$. For the nonlinear update, \eqref{eq:num-nonlinear-fit} and $H_{n,ii}+\rho\ge\rho$ yield
\[
|\alpha_{n,i}(x)| \le\frac{|b_{n,i}-(H_nx)_i|}{H_{n,ii}+\rho} \le C_R(1+\|x\|).
\]
The nonlinear update is globally Lipschitz on $\Omega_{n,R}$ by implicit differentiation, since
\[
\nabla_x\alpha_{n,i}(x) =-\frac{H_{n,i\cdot}} {H_{n,ii}+\rho+3\gamma\alpha_{n,i}(x)^2}.
\]
On the fixed time interval $[0,T]$, discrete Gronwall bounds fine, coarse, and a fixed number of corrected paths in a deterministic ball on $\Omega_{n,R}$, uniformly over the numerical grids. The first two state derivatives of the learner are uniformly bounded on such balls. Applying a smooth cutoff outside a larger ball produces globally $C_b^2$ coefficients agreeing with the original coefficients on the computed paths. Since the learner index takes finitely many values, its centered randomization can be represented with bounded, standardized categorical innovations. The population fields admit the corresponding local construction.

\subsection{Reference simulation for the second-order law}
\label{app:num-limit-reference}

The limit laws in Figure~\ref{fig:num-second-order-cdf} are generated from the population covariance $H=(0.5^{|r-s|})_{r,s=1}^4$ and the population training vector $b=1.6H\beta_v$, rather than any empirical $(H_n,b_n)$. For each learner, the limiting mean field $f$ and its first two derivatives follow directly from \eqref{eq:num-coordinate-learner} or \eqref{eq:num-nonlinear-fit}. Along the population mean flow $\theta$, each propagator has a joint centered Gaussian fluctuation and sensitivity pair $(\eta_j,\mathcal R_j)$ as defined in \eqref{eq:propagator-gaussian-increment}--%
\eqref{eq:propagator-random-sensitivity}. Its covariance is obtained by integrating the covariance of the transported learner fluctuations and their starting-state derivatives over the subinterval. These pairs are sampled independently across the $6$ subintervals, while the two components of each pair keep their joint covariance. The coarse Euler bias and its derivative are obtained by integrating \eqref{eq:propagator-euler-bias} along $\theta$. The propagator variables are inserted into the recursion \eqref{eq:general-limit-recursion} with $c=1$, which gives $50\,000$ independent draws of $[Z^{1,1}(T)]_1$ for each learner.

The covariance matrices are nonnegative up to integration roundoff, and the sample mean of the draws agrees, within Monte Carlo error, with the mean obtained by propagating the Euler bias through the recursion. The finite-$n$ values are the signed first-coordinate differences between the coupled Parareal and sequential fine endpoints, divided by $s_n^2$, from the $1000$ replications of each sample size.

\subsection{Timing details}\label{app:num-timing}

Table~\ref{tab:num-timing-processors} gives the elapsed times for all processor counts. The cost-model column uses unit costs measured on one subinterval: $c_F$ for a fine update, $c_{\rm tr}$ for a fine update that also stores the state, $c_V$ for the validation of one fine step, and $c_G$ for a coarse Euler step. It predicts the speed-up
\[
\frac{N(c_F+c_V)} {(K+1)MJc_G+(M/P)L\{Kc_F+c_{\rm tr}+c_V\}},
\]
which extends \eqref{eq:cost} by the validation costs and ignores communication. The last column is a second reference, which computes and stores the sequential path and then distributes only its validation over the $P$ processes. Storing the sequential path already costs more than the fused sequential run, so this reference is slower than the sequential one for both learners.

\begin{table}[htbp]
\centering
\small
\setlength{\tabcolsep}{4pt}
\caption{Parareal elapsed times (medians in milliseconds) and speed-ups relative to the sequential reference, for $P$ processes and the $M=6$ subintervals of Section~\ref{subsec:num-performance}. All configurations satisfy $E_n^K\le0.1\delta_n$ and $D_n^K\le0.1h_n^2$.} \label{tab:num-timing-processors}
\begin{tabular}{llrrrrr}
\toprule
Learner & $P$ & $K$ & \shortstack{Runtime\\(ms)} & \shortstack{Measured\\speedup} & \shortstack{Cost-model\\speedup} & \shortstack{Parallel\\validation only} \\
\midrule
Coordinate & 1 & 1 & 9.00 & 0.40 & 0.54 & 0.51 \\
&  & 2 & 10.79 & 0.33 & 0.46 &  \\
&  & 3 & 12.49 & 0.29 & 0.39 &  \\
\addlinespace[2pt]
& 2 & 1 & 4.97 & 0.72 & 1.08 & 0.63 \\
&  & 2 & 5.93 & 0.60 & 0.90 &  \\
&  & 3 & 7.00 & 0.51 & 0.78 &  \\
\addlinespace[2pt]
& 3 & 1 & 3.61 & 0.99 & 1.61 & 0.69 \\
&  & 2 & 4.35 & 0.82 & 1.34 &  \\
&  & 3 & 5.18 & 0.69 & 1.15 &  \\
\addlinespace[2pt]
& 6 & 1 & 2.51 & 1.42 & 3.15 & 0.71 \\
&  & 2 & 3.25 & 1.10 & 2.62 &  \\
&  & 3 & 3.71 & 0.96 & 2.24 &  \\
\midrule
Nonlinear & 1 & 1 & 38.30 & 0.45 & 0.53 & 0.85 \\
&  & 2 & 55.53 & 0.31 & 0.37 &  \\
&  & 3 & 72.72 & 0.24 & 0.28 &  \\
\addlinespace[2pt]
& 2 & 1 & 20.25 & 0.86 & 1.05 & 0.91 \\
&  & 2 & 29.25 & 0.59 & 0.72 &  \\
&  & 3 & 38.28 & 0.45 & 0.55 &  \\
\addlinespace[2pt]
& 3 & 1 & 14.07 & 1.24 & 1.57 & 0.92 \\
&  & 2 & 20.29 & 0.86 & 1.08 &  \\
&  & 3 & 26.57 & 0.65 & 0.82 &  \\
\addlinespace[2pt]
& 6 & 1 & 8.01 & 2.17 & 3.06 & 0.95 \\
&  & 2 & 11.63 & 1.50 & 2.10 &  \\
&  & 3 & 15.14 & 1.15 & 1.60 &  \\
\bottomrule
\end{tabular}
\end{table}

With one process, Parareal performs $K+1$ fine passes and is slower than the sequential run. The measured speed-ups grow with $P$, almost proportionally for the nonlinear learner, and stay below the cost model. The difference is largest for the coordinate learner, whose runs take a few milliseconds, so that the fixed cost of exchanging states between processes is not negligible.

\bibliographystyle{plainnat}
\bibliography{references}
\end{document}